\documentclass[11pt]{article}

\usepackage{epsf,latexsym,amssymb,graphics}
\usepackage[mathscr]{eucal}

\newtheorem{theorem}{Theorem}
\newtheorem{itlemma}{Lemma}[section]
\newtheorem{itproposition}[itlemma]{Proposition}
\newtheorem{itcorollary}[itlemma]{Corollary}
\newtheorem{itremark}[itlemma]{Remark}
\newtheorem{itremarks}[itlemma]{Remarks}
\newtheorem{itdefinition}[itlemma]{Definition}
\newtheorem{itexample}[itlemma]{Example}

\newenvironment{corollary}{\begin{itcorollary}\rm}{\end{itcorollary}} 
\newenvironment{definition}{\begin{itdefinition}\rm}{\end{itdefinition}}
\newenvironment{example}{\begin{itexample}\rm}{\end{itexample}}
\newenvironment{lemma}{\begin{itlemma}\rm}{\end{itlemma}} 
\newenvironment{proof}{\noindent {\em Proof}.\ \ }{\hspace*{\fill}$\halmos$\medskip}
\newenvironment{Proof}{\noindent }{\hspace*{\fill}$\halmos$\medskip}
\newenvironment{proposition}{\begin{itproposition}\rm}{\end{itproposition}}
\newenvironment{remark}{\begin{itremark}\rm}{\end{itremark}} 
\newenvironment{remarks}{\begin{itremarks} \rm}{\end{itremarks}}
\newenvironment{theo}{\begin{theorem} \rm}{\end{theorem}} 

\newcommand{\be}[1]{\begin{equation}\label{#1}}
\newcommand{\ee}{\end{equation}}

\newcommand{\bcl}{\begin{claim}}
\newcommand{\ecl}{\end{claim}}
\newcommand{\bc}[1]{\begin{corollary}\label{#1}}
\newcommand{\ec}{\end{corollary}}
\newcommand{\bd}[1]{\begin{definition}\label{#1}}
\newcommand{\ed}{\end{definition}}
\newcommand{\beqn}[1]{\begin{eqnarray}\label{#1}}
\newcommand{\eeqn}{\end{eqnarray}}
\newcommand{\beq}{\begin{eqnarray*}}
\newcommand{\eeq}{\end{eqnarray*}}
\newcommand{\bex}[1]{\begin{example}\label{#1}}
\newcommand{\eex}{\end{example}}
\newcommand{\bit}{\begin{itemize}}
\newcommand{\eit}{\end{itemize}}
\newcommand{\benu}{\begin{enumerate}}
\newcommand{\eenu}{\end{enumerate}}
\newcommand{\bfa}[1]{\begin{fact}\label{#1}}
\newcommand{\efa}{\end{fact}}
\newcommand{\bl}[1]{\begin{lemma}\label{#1}}
\newcommand{\el}{\end{lemma}}
\newcommand{\bpr}{\begin{proof}}
\newcommand{\epr}{\end{proof}}
\newcommand{\bpro}{\begin{Proof}}
\newcommand{\epro}{\end{Proof}}
\newcommand{\bp}[1]{\begin{proposition}\label{#1}}
\newcommand{\ep}{\end{proposition}}
\newcommand{\rf}[1]{~(\ref{#1})}
\newcommand{\br}[1]{\begin{remark}\label{#1}}
\newcommand{\er}{\end{remark}}
\newcommand{\brs}[1]{\begin{remarks}\label{#1}}
\newcommand{\ers}{\end{remarks}}
\newcommand{\bt}[1]{\begin{theo}\label{#1}}
\newcommand{\et}{\end{theo}}

\newcommand{\comment}[1]{}

\newcommand{\halmos}{\rule{1ex}{1.4ex}}
\newcommand{\qedb}{\hfill \halmos}

\newcommand{\abs}[1]{|{#1}|}
\newcommand{\dnorm}[1]{\|{#1}\|}

\newcommand{\C}{{\mathbb C}}
\newcommand{\N}{{\mathbb N}}  
\newcommand{\R}{{\mathbb R}}  

\newcommand{\U}{{\mathcal U}}

\newcommand{\flr}[1]{\lfloor #1 \rfloor}
\newcommand{\clg}[1]{\lceil #1 \rceil}
\newcommand{\eps}{{\varepsilon}}
\newcommand{\ap}{{\mathcal A}}

\newcommand{\cde}{\mathcal C}
\newcommand{\cdek}{{\bf C}_{\mbox{\tiny $K,\ell$}}}
\newcommand{\cmax}{{\bf C}_{\mbox{\tiny max}}} 
\newcommand{\copt}{{\bf C}_*}
\newcommand{\fbk}{\mbox{\tiny fb}}
\newcommand{\tot}{\mbox{\tiny tot}}

\begin{document}

\title{Optimal Length and Signal Amplification in\\
Weakly Activated Signal Transduction Cascades}

\author{Madalena Chaves$^{1,2}$,
Eduardo D.\ Sontag$^1$
and Robert J.\ Dinerstein$^2$\\
{\small $^1$ Department of Mathematics, Rutgers University, 
New Brunswick, NJ 08903}\\
{\small $^2$ Lead Generation Informatics, Aventis, Bridgewater, NJ 08807} 
}

\date{}

\maketitle

\begin{abstract}
Weakly activated signaling cascades can be modeled as linear
systems. The input-to-output transfer function and the internal 
gain of a linear system, provide natural  measures for the 
propagation of the input signal down the cascade and for the 
characterization of the final outcome.
The most efficient design of a cascade for generating sharp signals,
is obtained by choosing all the off rates equal, and 
a ``universal'' finite optimal length.
\end{abstract}




\section{Introduction}
\label{sec-introduction}
Protein kinase cascades are major functional modules used by cells to
translate signals generated by receptor activation into diverse biochemical
and physiological responses~\cite{pl-ejb}.  Highly conserved
throughout evolution and across species, the kinase cascade motif
participates in the control of many processes, including cell cycle
regulation, gene expression, cellular metabolism, stress responses, and T
cell activation.  For this reason, control of kinase cascades by therapeutic
intervention has become an attractive area for drug discovery, particularly
in the areas of cancer and inflammation~\cite{cohen-nrdd,aventis}.

Some four mitogen-activated protein kinase (MAPK) signaling cascades have
been found in yeast~\cite{pg-science} and at least a dozen
MAPK cascades have been identified in mammalian cells~\cite{karan}.  
The intensive study of MAPK
pathways has prompted efforts to characterize these systems theoretically
(see, inter alia,~\cite{asthagiri,bhalla-chaos,bhalla-science,
kholo-jtb,hnr2002,huang-ferrel,kholo-jeb,kholo-febs,levchenko}).  
In this paper, we will
utilize concepts and methods from the theory of linear control systems
to characterize kinase signaling cascades, and in particular the MAPK
pathway, in order to understand how the number of kinases in a cascade
and their individual enzymatic activities can affect the pathway in its role as a
signal transducing module.

Let $R$ denote the input signal, 
$\tilde X_i$  the inactive (nonphosphorylated) form of kinase
$i$ and $X_i$ the active (phosphorylated) form of kinase $i$. 
The rate constant (or ``on'' rate) for the $i$-th kinase phosphorylation  
will be denoted by $\tilde\alpha_i$, and
the dephosphorylation rate constants (or ``off'' rate) will be denoted
$\beta_i$.
The input signal $R$ might represent, for example, the concentration of activated
receptors, and the dynamics of the signal transduction pathway may be
modeled as follows (see~\cite{hnr2002}):
\beqn{mapk-cascade0}
   \frac{d X_1}{dt} = \tilde\alpha_1R\tilde X_1 -\beta_1X_1,
   \ \ \ \ 
   \frac{d X_i}{dt} = \tilde\alpha_iX_{i-1}\tilde X_i -\beta_iX_i,
   \ \ i=2,\ldots,n.
\eeqn
Assuming that the total amount of kinase $i$ remains constant, that is, 
$
   \tilde X_i+X_i=X_{\tot,i}
$
the differential equations\rf{mapk-cascade0} can be rewritten as
\beq
   \frac{d X_1}{dt} &=& \alpha_1 R
   \left(1-\frac{X_1}{X_{\tot,1}}\right) -\beta_1X_1  
\eeq
and
\beqn{mapk-cascade}  
   \frac{d X_i}{dt} &=& \alpha_i X_{i-1}
   \left(1-\frac{X_i}{X_{\tot,i}}\right) -\beta_iX_i,  
   \ \ i=2,\ldots,n.
\eeqn
where $\alpha_i=\tilde\alpha_i X_{\tot,i}$.
Throughout this paper we will focus on the case of
{\it weakly activated pathways}, by which we mean a 
low level of kinase phosphorylation, that is
\beqn{eq-weak-act}
   X_i\ll X_{\tot,i}\ \ \ \Rightarrow\ \ \ 1-\frac{X_i}{X_{\tot,i}}\approx 1.
\eeqn 
In this case the equations\rf{mapk-cascade} are simplified to a linear
system of the form: 
\beqn{eq-linear}
   \frac{d X_1}{dt} = \alpha_1 R -\beta_1X_1,   
\ \ \
  \frac{d X_i}{dt} = \alpha_i X_{i-1} -\beta_iX_i,\ \ \ i=2,\ldots n.
\eeqn
In Section~\ref{sec-tf} we will describe how to compute the transfer
function and internal gain for this system and then in 
Section~\ref{sec-definitions} we will define a set of measures for
the output signal, which closely follow those discussed in~\cite{hnr2002}.
In Section~\ref{sec-design} we prove that the most efficient cascade
design, for generating sharp signals, has equal on rates and 
a finite length depending only on the cascade's internal gain.
In Section~\ref{sec-feedback} positive feedback from the last
activated kinase to the first is added to the cascade, and the 
optimal design is re-examined in this new context.
Finally, in Sections~\ref{sec-delay} and~\ref{sec-stability} we briefly discuss 
the effect of delays along the cascade and how to check the cascade's stability
to random small perturbations.

\section{The input-to-output transfer function}
\label{sec-tf}

We will consider the signaling cascade\rf{eq-linear} as a system with
an {\it input $R$}, and an {\it output} which will be some function of the 
concentration of the last kinase $X_n$.  
Specifically, the output will be the ``effective'' integral of $X_n$, or
in other words, the cascade will be extended one more step to
include a ``leaky'' integrator:
\beq
   \frac{d X_{n+1}}{dt}=X_{n}-\ell X_{n+1},
\eeq
where the output is $X_{n+1}$.
The variable $X_{n+1}$ expresses the effective concentration of the last kinase
(minus losses due to degradation or inactivation of $X_n$, for instance). 
Note that the case $\ell=0$ recovers $X_{n+1}=\int^t X_n(t')\; dt'$.
\begin{figure}[h,t]
\epsfysize=8cm
\epsfxsize=7cm
\centerline{\epsffile{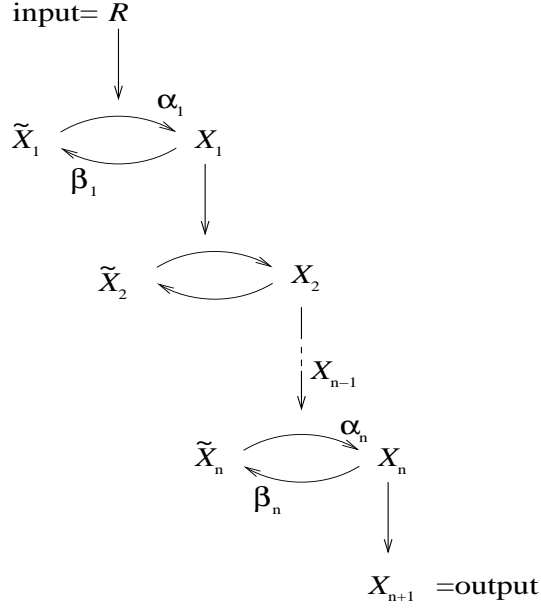}}
\caption{A model of a MAPK cascade. }
\label{fig-mapk}
\end{figure}

The model for a weakly activated signal transduction cascade may then
be written in the more compact form, 
\beqn{eq-in-out-sys}
   \frac{d X}{dt}(t) = A X(t) + B R(t), 
   \ \ \ Y(t)= C X(t),
\eeqn
where $X=(X_1,X_2,\ldots,X_n,X_{n+1})'$ is a column vector whose elements are
always nonzero, and $A\in\R^{(n+1)\times(n+1)}$,
$B\in\R^{(n+1)\times1}$ and $C\in\R^{1\times(n+1)}$ are the matrices
\beqn{eq-AB}
    A=\pmatrix{-\beta_1 & 0 & 0 & \cdots & 0 & 0 & 0\cr
               \alpha_2 & -\beta_2 & 0 & \cdots & 0 & 0 & 0\cr
               0 & \alpha_3 & -\beta_3 & \cdots & 0 & 0 & 0\cr
               0 & 0 & \alpha_4 & \cdots & 0 & 0 & 0\cr
               \vdots &  & & \ddots &  & &\vdots \cr
               0 & 0 & 0 & \cdots &\alpha_n & -\beta_n  & 0\cr
               0 & 0 & 0 & \cdots & 0  & 1  & -\ell
               },\ \ \
   B=\pmatrix{\alpha_1\cr 0\cr 0\cr 0\cr \vdots \cr 0 \cr 0}, 
\eeqn
and  
\beqn{eq-C}
   C=\pmatrix{0 & \cdots & 0 & 1}.
\eeqn

It is well known (see~\cite{control,mct}, or any other book on control systems)
that, for a system such as\rf{eq-in-out-sys},  the output can be 
computed directly as the convolution between the input signal $R$ and the
impulse response of the system. The impulse response of the system is 
the output corresponding to a single input pulse. 
If we let $G$ denote the impulse response (and
assuming that the system starts at rest, with initial condition $X(0)=0$), then
\beq
    Y(t) = (G*R)\ (t).
\eeq
The impulse response, $G$,  characterizes the action that the internal
structure of the system will have on any input, such as the filtering
of certain frequency components, and the amplification or dampening of 
the signal.
Biological inputs may take many different forms, such as single pulses,
slowly decaying signals, constant stimuli applied for a
certain time interval, or oscillatory signals. 
Thus, it is appropriate to have a model in which the output signal is
obtained as a convolution of the input $R$ (which may take many forms)
and the transfer function $G$ (which depends only on the intrinsic
kinase activity parameters, and needs to be computed only once).

A very convenient way to analyze system\rf{eq-in-out-sys} is to
convert it to the {\it frequency domain}, by application of the
Laplace transform operator.
The Laplace transform of the impulse response is called the {\it
transfer function} of the system, and it provides a simple linear
relationship between the Laplace transforms of the input and the
output, as well as also providing a measure of amplification/dampening 
of the input signal.
The transfer function is given by a simple formula in terms of the matrices 
$A$, $B$ and $C$ as summarized in Appendix~\ref{sec-fd}. For this
cascade system, we will carry out the Laplace transforms in detail so as to gain 
some insight into the internal structure of the system.

The Laplace transform of $X$ will be denoted by $\hat X$, and is
defined as
\beq
   \hat X_i(s) = \int_{-\infty}^{\infty} e^{-st} X_i(t) dt, \
   \ \mbox{ and }\ \  
   \hat R(s) = \int_{-\infty}^{\infty} e^{-st} R(t) dt,
\eeq
where $s$ is a complex number $s=s_{re}+\jmath\omega$ ($\jmath$ is the
imaginary number $\sqrt{-1}$) and takes values in an appropriate
region of convergence.

Applying the Laplace transform operator to both sides of
equations\rf{eq-linear},  assuming that $X(0)=0$, and recalling the
properties of the Laplace transform (see Appendix~\ref{sec-fd}), we have:
\beq
   s \hat X_1(s) &=& \alpha_1\hat R(s)-\beta_1\hat X_1(s)\\
   s \hat X_i(s) &=& \alpha_i\hat X_{i-1}(s)-\beta_i\hat X_i(s),
   \ \ \ i=2,\ldots,n \\
   s \hat X_{n+1}(s) &=& \hat X_{n}(s)-\ell\hat X_{n+1}(s),
\eeq
which yields
\beq
   \hat X_1(s) = \frac{\alpha_1}{s+\beta_1}\;\hat R(s),\ \ \ 
   \hat X_i(s) = \frac{\alpha_i}{s+\beta_i}\;\hat X_{i-1}(s),
   \ \ \ i=2,\ldots,n,
\eeq
and 
\beq
    \hat X_{n+1}(s) = \frac{1}{s+\ell}\;\hat X_n(s).
\eeq
In this way, we may view the cascade as a sequence of $n$ steps, the output of
the step $i-1$ becoming the input to step $i$. 
\begin{figure}[h,t]
\epsfysize=0.8cm
\epsfxsize=12cm
\centerline{\epsffile{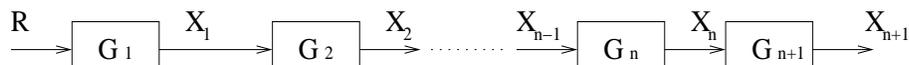}}
\caption{Transfer functions at each step. }
\label{fig-cascade}
\end{figure}

For each {\it single step} in the cascade, the input is $\hat X_{i-1}$
and the output is $\hat X_{i}$, and they are related by a multiplicative
factor, which is in fact the {\it transfer function for the step $i$}:
\beq
   \hat G_i(s)=\frac{\alpha_i}{s+\beta_i},\ \ \ i=2,\ldots,n.
\eeq
For the whole cascade, the input is $R$ and the output is $\hat X_{n+1}$, and
it is easy to see that the transfer function for the total system is
the product of the transfer functions at each step:
\beqn{eq-Ghat}
   \hat G(s)=\hat G_1(s)\cdots \hat G_{n+1}(s)=
   \frac{1}{s+\ell}\frac{\alpha_1\cdots\alpha_n}{(s+\beta_1)\cdots(s+\beta_n)}.
\eeqn
Therefore
\beq
   \hat Y(s)=\hat X_{n+1}(s) =
   \frac{1}{s+\ell}\frac{\alpha_1\cdots\alpha_n}{(s+\beta_1)\cdots(s+\beta_n)}
   \hat R(s),
\eeq
and the actual output may now be obtained by the inverse Laplace
transform.
Alternatively, even without knowing the exact form of the output, that
is, the function $Y(t)$, 
it is still possible to further characterize the properties 
of the system, through the {\it 2-norm} of the functions $\hat Y$ and 
$\hat R$. Define the {\it 2-norm} of the function $Y$ and  the 
{\it 2-norm} of the Laplace transform $\hat Y$ by
\beq
   \dnorm{Y}_2 := \left[\int_{-\infty}^{+\infty} 
                      \abs{Y(t)}^2\ dt\right]^{\frac{1}{2}},
   \ \ \mbox{ and }\ \ 
   \dnorm{\hat Y}_2 := \left[\frac{1}{2\pi}\int_{-\infty}^{+\infty} 
                      \abs{\hat Y(\jmath\omega)}^2 d\omega\right]^{\frac{1}{2}},
\eeq
and similar expressions hold for $R$ and $\hat R$.
(Note: from now on we will assume that the signals are defined only
for positive times, that is, $Y(t)=0$ for $t<0$).
The 2-norm $\dnorm{Y}_2$ provides a measure of the {\it strength} of the signals 
(in analogy to the energy of a mechanical system).
Indeed, these norms provide a very convenient way to relate the input and 
output because, from Parseval's Theorem, 
{\it the 2-norm of a function, equals the 2-norm of its Laplace transform},
and therefore
\beq
  \dnorm{\hat Y}_2=\dnorm{Y}_2,\ \ \ \ \dnorm{\hat R}_2=\dnorm{R}_2,
\eeq
without the need to compute inverse transforms (a very helpful
fact, since in general the inverse transforms may not be simple to compute).

Another useful measure is the {\it infinity norm} of the transfer
function, that selects the least upper bound of the absolute value of
$\hat G$,
\beq
   \dnorm{\hat G}_{\infty}:=\sup_{\omega\in\R} \abs{\hat G(\jmath\omega)}.
\eeq
As shown in Appendix~\ref{sec-fd}, a very useful estimate for characterizing
the relative strength of the input and output signals is:
\beqn{eq-L2-ineq}
    \dnorm{Y}_2 \leq \dnorm{\hat G}_{\infty} \dnorm{R}_2,
\eeqn
where it is immediately apparent that the infinity norm of the 
transfer function gives an upper bound for the amplification of the
input signal throughout the cascade. 
Moreover, the infinity norm $\dnorm{\hat G}_{\infty}$ is in fact the
smallest number that satisfies\rf{eq-L2-ineq}, for all input/output
pairs (that is, pairs $(R,Y)$, where $Y$ is the output corresponding to
the input $R$).

To compute the infinity norm of the transfer function for the whole
cascade, note that
\beq
   \abs{\hat G_i(\jmath\omega)}^2 = 
   \frac{\abs{\alpha_i}^2}
        {\abs{\jmath\omega+\beta_i}^2}
   \equiv \frac{\alpha_i^2}{\omega^2+\beta_i^2}
   \leq   \frac{\alpha_i^2}{\beta_i^2}, 
    \ \ \ \mbox{for all }\omega\in(-\infty,\infty),   
\eeq 
and the equality holds for $\omega=0$. Therefore
\beqn{eq-Ginf}
   \dnorm{\hat G}_{\infty} = \frac{1}{\ell}\frac{\alpha_1\cdots\alpha_n}
                                  {\beta_1\cdots\beta_n}.
\eeqn
A necessary condition for amplification of the signal to occur is that
$\dnorm{\hat G}_{\infty}>1$. Moreover, since $\ell$ is essentially an
independent parameter, introduced 
for the purpose of defining a reasonable measure of the output, we 
can say that amplification of the input signal occurs only if 
\beqn{eq-amplf}
    \alpha_1\cdots\alpha_n > \beta_1\cdots\beta_n.
\eeqn
Recall that $\alpha_i\equiv\tilde\alpha_i X_{\tot,i}$, where $X_{\tot,i}$ is the
total concentration of the $i$th kinase and $\tilde\alpha_i$ is the 
(true) rate of phosphorylation. Therefore, we still expect that
$\tilde\alpha_i<\beta_i$, $i=1,\ldots,n$, as should be the case for a weakly 
activated pathway.
 
The norm $\dnorm{\hat G}_{\infty}$ is often called the 
{\it internal gain of the system} which,
through expression\rf{eq-L2-ineq}, provides a useful and easy way to compute
the input-to-output strength relation.
For example, if a MAPK cascade has a ``5-fold amplification'', then 
its internal gain is $\dnorm{\hat G}_{\infty}=5$.

Note that, in the case where $\ell=0$, the internal gain $\dnorm{\hat G}_{\infty}$ is
infinite --- meaning that, in at least one step (
$X_n\to X_{n+1}$ ) there is no degradation term.
Then the estimate\rf{eq-L2-ineq} contains no useful information.
However, for $\ell=0$, we have
\beq
   Y(t) = X_{n+1}(t)=\int_{0}^{t}\ X_n(t')\ dt',
\eeq
and we also have an estimate for the ``strength'' of the signal $X_n$, since
\beq
   \dnorm{X_n}_2 \leq \frac{\alpha_1\cdots\alpha_n}
                           {\beta_1\cdots\beta_n} \dnorm{R}_2.
\eeq

\section{Signaling time, signal duration and signal amplitude}
\label{sec-definitions}

Some basic quantities which serve to characterize a signal
transduction system are: the overall amplification from the input 
to the ouput; the duration of the output signal; and the time it 
takes the input signal to traverse the cascade.
There are several possible definitions and estimates of these
quantities: here we extend the definitions given by~\cite{hnr2002},
embedding them in the context of frequency-domain analysis, and 
generalizing them to arbitrary inputs.

To be concise, let us identify the cascade\rf{eq-in-out-sys} by its parameters, 
and associate with it the following $(2n+1)$-tuple: 
\beq
    \cde := (n,\alpha_1,\ldots,\alpha_n,\beta_1,\ldots,\beta_n),
\eeq
where it is assumed that $n\in\N$ and $\alpha_i,\beta_i$ are positive real 
numbers, for $i=1,\ldots,n$.
We will also introduce the notation $\U$ for denoting the set of inputs. 

\bd{def-sigtime}
For system\rf{eq-in-out-sys}, with parameters $\cde$ and a leak factor $\ell>0$,
for each input $R$,
the {\it signaling time}, $\tau$, and the {\it output signal
duration}, $\sigma$, are given by
\beq
   \tau(\cde,\ell,R):=-\left.\frac{d\; \ln\hat Y}{ds}(s)\right\rfloor_{s=0},
\ \ \ \ \ 
   \sigma(\cde,\ell,R):=\sqrt{\left.\frac{d^2\; \ln\hat Y}{ds^2}(s)\right\rfloor_{s=0}
   }\ .
\eeq
The {\it signaling time to step i} 
and the {\it signal duration at step i}, $i\leq n$ are given by:
\beq
   \tau_i(\cde,R):=-\left.\frac{d\; \ln\hat X_i}{ds}(s)\right\rfloor_{s=0},
\ \ \ \ \ 
   \sigma_i(\cde,R):=\sqrt{\left.\frac{d^2\; \ln\hat X_i}{ds^2}(s)\right\rfloor_{s=0}
   }\ .
\eeq
\ed
To understand the significance of these definitions, recall the properties of 
the Laplace transform and compute (with $Y(t)=0$ for $t\leq0$):
\beq
  \hat Y(0) = \int_0^\infty Y(t) dt, \ \ \
  \frac{d\hat Y}{ds}(0) = -\int_0^\infty t Y(t) dt,\ \ \
  \frac{d^2\hat Y}{ds^2}(0) = \int_0^\infty t^2 Y(t) dt 
\eeq
and thus we recover expressions (4) and (5) of reference~\cite{hnr2002}
\beq 
   \tau = \frac{\int_0^\infty t Y(t) dt}
                 {\int_0^\infty  Y(t) dt}, \ \ \ \ 
   \sigma^2 = \frac{\int_0^\infty t^2 Y(t) dt}{\int_0^\infty Y(t) dt}
    -\left(\frac{\int_0^\infty t Y(t) dt}{\int_0^\infty Y(t) dt}\right)^2,
\eeq
where $\tau$ can be regarded as the expected value (of the time to
traverse the pathway), and $\sigma$ as the corresponding variance.

An estimate of the amplitude of the output signal, as given in equation (6) of 
reference~\cite{hnr2002}, is the value $S$, such that 
$S\times 2\sigma=\int_0^\infty Y(t) dt$.
Again we propose a more generalized notion, suggested by the
input-to-output estimate\rf{eq-L2-ineq}, that takes advantage of
the easily computed internal gain  of the system,
and also incorporates the strength of the signal. 

\bd{def-amplitude}
For system\rf{eq-in-out-sys}, with parameters $\cde$ and a leak factor $\ell>0$,
for each input $R$, the {\it signal amplitude} is given by
\beqn{eq-amplitude}
   \ap(\cde,\ell,R)\ :=\ \frac{\dnorm{\hat G_{\cde}}_\infty}{\sigma(\cde,\ell,R)}\ \dnorm{R}_2,
\eeqn
where $\hat G_{\cde}$ is the transfer function\rf{eq-Ghat}.
\ed

$\ap$ may also be regarded as the amplitude of a constant signal of duration $\sigma$, 
but Definition~\ref{def-amplitude} differs from the definition of amplitude given 
in~\cite{hnr2002} in essentially three points:
\bit
\item[1.] the meaningful quantity for measuring the amplitude is not the integral 
$\int Y(t)\;dt$ (which computes the area under the curve $Y(t)$), but rather 
the 2-norm  $\sqrt{\int \abs{Y(t)}^2 dt}$, which computes the strength 
of the signal;

\item[2.] the amplitude is proportional to the product of 
the internal gain of the system, and the 2-norm of the input. 
This simplifies calculations since, for each cascade,
the $\dnorm{\hat G}$ is computed only once and $\dnorm{R}_2$ is computed for each 
input signal;

\begin{figure}[h,t]
\epsfysize=5cm
\epsfxsize=6cm
\centerline{\epsffile{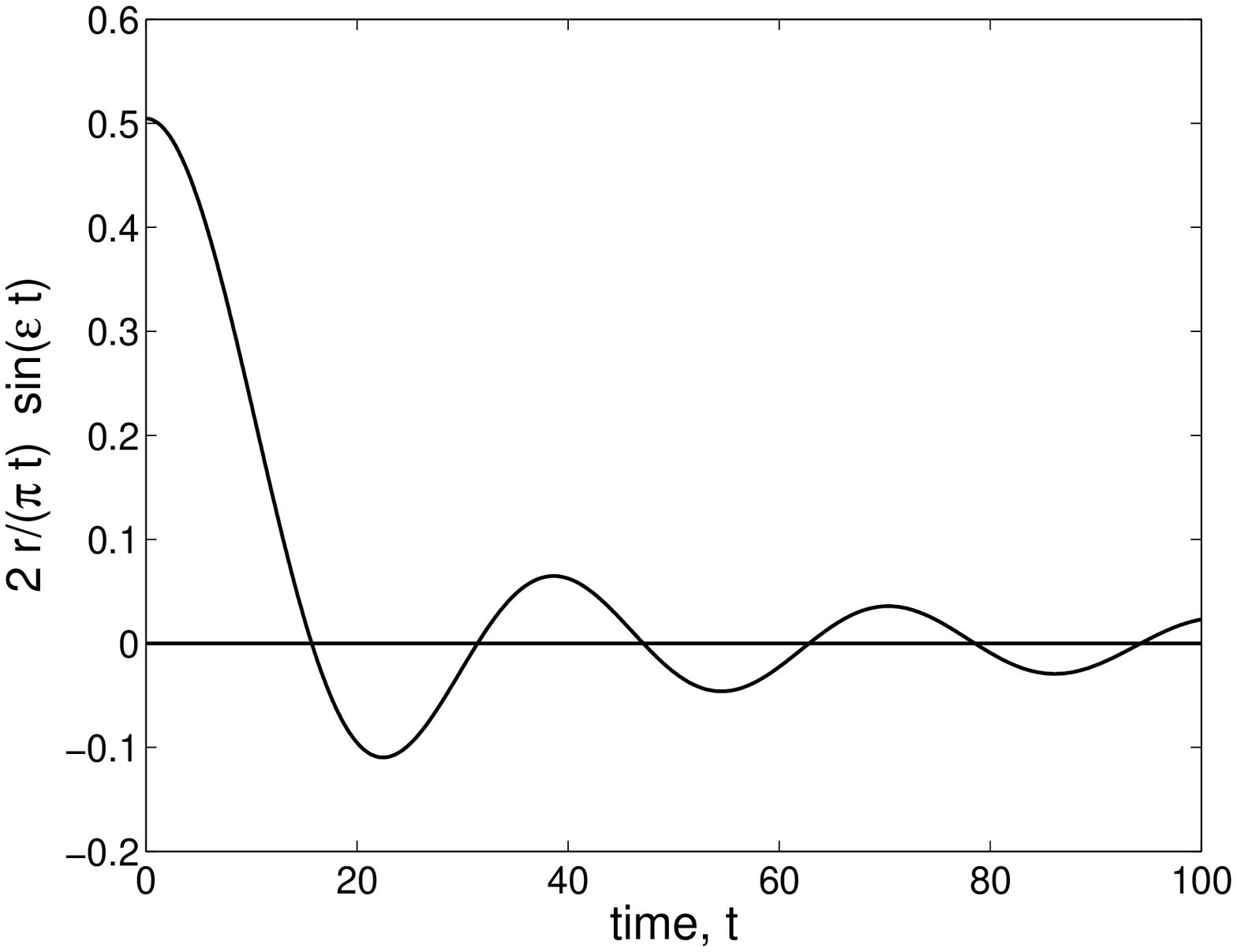}}
\caption{An input that satisfies $\dnorm{R}_2=1$ and 
    $\dnorm{Y}_2\approx\dnorm{\hat G}_\infty$, with $\eps=0.2$.  }
\label{fig-sync}
\end{figure}

\item[3.] the product $\dnorm{\hat G}_\infty\dnorm{R}_2$ is used as an estimate for 
$\dnorm{Y}_2$, but we know (see Appendix) that $\dnorm{\hat G}_\infty$ is 
the least factor that satisfies the inequality $\dnorm{Y}_2\leq \kappa\dnorm{R}_2$.
In fact,~\cite{control2} shows how to construct examples of inputs for which the 
equality is approximated. For instance, for any $\eps>0$,
the input depicted in Figure~\ref{fig-sync}:
\beqn{eq-sync}
     R(t)=2\frac{r}{\pi t}\sin{\eps t}, \ \ \mbox{ with }\ \ r=\sqrt{\pi/\eps},
    \ \ \mbox{ for }\ \ t\geq0,
\eeqn
has unit norm, i.e., $\dnorm{R}_2=1$, and satisfies 
$\dnorm{Y}_2\approx\dnorm{\hat G}_\infty$, for $\eps$ small enough,
as shown in the Appendix.
\eit

We remark that these definitions are valid not only for the 
special case when $A$, $B$ and $C$ are of the form specified
in equations\rf{eq-AB},\rf{eq-C}, 
but in fact they are valid for any linear system of
the form\rf{eq-in-out-sys}. For instance, in
Section~\ref{sec-feedback}, we compute these quantities for 
the case when there is positive feedback from the last to the first kinase.
We next explicitly compute these quantities for the special case 
when $A$, $B$ and $C$ are of the form\rf{eq-AB} and\rf{eq-C}, and $\ell=0$:
\beqn{eqc-tau}
    \tau(\cde,\ell,R)=\frac{1}{\ell}
         +\sum_{i=1}^{n}\frac{1}{\beta_i} 
         + \left.\frac{d\ \ln\hat R}{ds}\right\rfloor_{s=0}
\eeqn
\beqn{eqc-sigma}
    \sigma(\cde,\ell,R)=
    \sqrt{\frac{1}{\ell^2}+\sum_{i=1}^{n}\frac{1}{\beta_i^2} + q(R)},
    \ \ \mbox{ where }\ \ 
    q(R)=\left.\frac{d^2\ \ln\hat R}{ds^2}\right\rfloor_{s=0}
\eeqn
\beqn{eqc-amplitude}
   \ap(\cde,\ell,R)
       =\frac{1}{\sqrt{\frac{1}{\ell^2}+\sum_{i=1}^{n}\frac{1}{\beta_i^2} + q(R)}}\
       \frac{\alpha_1\cdots\alpha_n}{\ell\beta_1\cdots\beta_n}\ \dnorm{R}_2\ .
\eeqn
In the case $\ell=0$, the quantities $\tau$, $\sigma$ and $\ap$ may be
computed for $Y\equiv X_n$. The expressions are very similar, except that all
the terms in $\ell$ vanish.

\bex{ex-inputs} 
A typical input is a decaying exponential $R(t)=R_0e^{-\lambda t}$,
with
\beq
  \dnorm{R}_2 = \frac{R_0}{2\lambda},\ \ \ 
  \hat R (s)=\frac{R_0}{s+\lambda}, \ \ \ 
  \frac{d\ \ln\hat R}{ds}\ (0)=-\frac{1}{\lambda},\ \ \ 
  q(R)=\frac{1}{\lambda^2}.
\eeq

A ``peak''-like input may be represented by $R(t)=R_0 t e^{-\lambda t}$,
with
\beq
  \dnorm{R}_2 = \frac{R_0}{4\lambda^3},\ \ \ 
  \hat R (s)=\frac{R_0}{(s+\lambda)^2}, \ \ \ 
  \frac{d\ \ln\hat R}{ds}\ (0)=-\frac{2}{\lambda},\ \ \ 
  q(R)=\frac{2}{\lambda^2}.
\eeq

For a constant signal, of magnitude $R_0$, applied for an interval of
time $T_0$, we have:
\beq
  \dnorm{R}_2 = R_0 \sqrt{T_0},\ \ \ 
  \hat R (s)=R_0\,\frac{1-e^{-sT_0}}{s}, \ \ \ 
  \frac{d\ \ln\hat R}{ds}\ (0)=-\frac{T_0}{2},\ \ \ 
  q(R)=\frac{T_0^2}{12}.
\eeq
\eex

\section{Cascade design optimization}
\label{sec-design}
From the analysis of the quantities $\tau$, $\sigma$ and $\ap$,
defined in Section~\ref{sec-definitions}, we can explore the
signaling efficiency of kinase cascades.
The definition of an ``efficient'' response may depend on the
particular biological context, but it typically involves the relationship
between the length of the cascade, the amplitude of the signal
and its duration.
A question posed in~\cite{hnr2002} is whether cascades can respond 
with sharp signals, i.e., simultaneously of short duration and high amplitude. 
Our model provides a definite answer to this question.

As we have seen, our linear model has a gain
that depends on the length of the cascade and the values of
the on/off rate constants, but doesn't depend on the input.
As a starting point, we may think of the family of cascades that
have the same value for the internal gain, say $K$, and examine their
length, the distribution of the ``on/off'' rates and signal 
amplitude and duration.
The problem we would like to study is then:
\bit
\item[(P)]
For each fixed internal gain, 
$\dnorm{\hat G}_\infty=K$, find the optimal combination of the on/off
rates and the length of the cascade that maximizes the signal
amplitude, $\ap$, for any input $R$.
\eit

To formulate this problem, first define the family of cascades that 
have the same internal gain $K$:
\beq
   \cdek:=\{ \cde=(n,\alpha_1,\ldots,\alpha_n,\beta_1,\ldots,\beta_n):\ 
                  \frac{\alpha_1\cdots\alpha_n}{\ell\beta_1\cdots\beta_n}=K \},
\eeq
For each input $R$, and each leak factor $\ell$, define the set of 
``optimal'' cascades, that is, those cascades which exhibit maximal 
signal amplitude:
\beq
   \cmax(\ell,R) := \{ \cde\in\cdek:\ \ap(\cde,\ell,R)\geq \ap(\cde',\ell,R),\ 
                      \mbox{ for all }\cde'\in\cdek \}.  
\eeq
Then define the function
\beq
    \sigma_0(n,\beta_1,\ldots,\beta_n)
       :=\sum_{i=1}^{n}\frac{1}{\beta_i^2} 
\eeq
and observe that it satisfies
\beq
    \sigma(\cde,\ell,R)=\sqrt{\frac{1}{\ell^2}
                 +\sigma_0(n,\beta_1,\ldots,\beta_n)+q(R)}.
\eeq
Finally, define the set of cascades that minimize $\sigma_0$ over the family $\cdek$:
\beq
    \copt(\ell,R) := \{ \cde\in\cdek:\ \sigma_0(n,\beta_1,\ldots,\beta_n)
                                  \leq \sigma_0(n',\beta_1',\ldots,\beta_n'),\ 
                      \mbox{ for all }\cde'\in\cdek  \}.
\eeq
Our first result states that in fact the sets $\copt(\ell,R)$ and 
$\cmax(\ell,R)$ are equal, or in other words, that an optimal cascade
will {\it simultaneously maximize the signal amplitude and minimize the
  signal duration}.

\bl{lm-Rindep}
In the notation defined above, $\cmax(\ell,R)=\copt(\ell,R)$, for all inputs
$R\in\U$ and leak factors $\ell>0$. 
\el

\bpr
Fix any $\ell>0$, and any $R\in\U$.
Recall the notation $\cde=(n,\alpha_1,\ldots,\alpha_n,\beta_1,\ldots,\beta_n)$. 
Given any $\cde,\cde'\in\cdek$, the following equivalences hold:
\beqn{equiv-sigma0}
   &&  \sigma_0(n,\beta_1,\ldots,\beta_n)
    \leq \sigma_0(n',\beta_1',\ldots,\beta_n'),  \nonumber \\ 
   \ \ & \Leftrightarrow & \ \    \nonumber
    \sqrt{ \frac{1}{\ell^2}+\sigma_0(n,\beta_1,\ldots,\beta_n)+q(R)} 
    \leq \sqrt{\frac{1}{\ell^2}+\sigma_0(n',\beta_1',\ldots,\beta_n')+q(R)},\\ 
    \ \ & \Leftrightarrow & \ \ 
  \sigma(\cde,\ell,R)\leq \sigma(\cde',\ell,R)
\eeqn
and also
\beqn{equiv-ap}
   \sigma(\cde,\ell,R)\leq \sigma(\cde',\ell,R) 
   \ \ & \Leftrightarrow & \ \  \nonumber
   \frac{K\ \dnorm{R}_2}{\sigma(\cde,\ell,R)}
   \geq \frac{K\ \dnorm{R}_2}{\sigma(\cde',\ell,R)} \\
    \ \ & \Leftrightarrow & \ \ 
   \ap(\cde,\ell,R)\geq\ap(\cde',\ell,R).
\eeqn
Therefore,\rf{equiv-sigma0} and\rf{equiv-ap} imply that, for any 
two cascades $\cde,\cde'\in\cdek$,
\beqn{equiv-ap-sig0}
    \sigma_0(n,\beta_1,\ldots,\beta_n)
    \leq \sigma_0(n',\beta_1',\ldots,\beta_n')\ \ 
   \ \ & \Leftrightarrow & \ \   
     \ap(\cde,\ell,R)\geq\ap(\cde',\ell,R).
\eeqn
To show that $\copt(\ell,R)$ is contained in $\cmax(\ell,R)$, 
pick any $\cde\in\copt(\ell,R)$. Then
\beq
  \sigma_0(n,\beta_1,\ldots,\beta_n)
    \leq \sigma_0(n',\beta_1',\ldots,\beta_n'),
 \mbox{ for all } \cde'\in\cdek.
\eeq
By\rf{equiv-ap-sig0}, this is equivalent to 
$\ap(\cde,\ell,R)\geq\ap(\cde',\ell,R)$, for all  $\cde'\in\cdek$,
and so  $\cde\in\cmax(\ell,R)$.

Conversely, we need to show that $\cmax(\ell,R)$ is contained in $\copt(\ell,R)$.
So, pick any $\cde\in\cmax(\ell,R)$. It satisfies:
\beq
   \ap(\cde,\ell,R)\geq \ap(\cde',\ell,R),\ \ \mbox{ for all }\cde'\in\cdek.
\eeq
Again by\rf{equiv-ap-sig0}, this is equivalent to 
$\sigma_0(n,\beta_1,\ldots,\beta_n)
  \leq \sigma_0(n',\beta_1',\ldots,\beta_n')$
for all $\cde'\in\cdek$.
We conclude that $\cde\in\copt(\ell,R)$, as we wanted to show.
\epr

An immediate conclusion from Lemma~\ref{lm-Rindep} is that, 
\beq
   \mbox{ maximize $\ap(\cde,\ell,R)$ over $\cdek$}
   \ \ \Leftrightarrow  \ \ 
   \mbox{ minimize $\sigma_0(n,\beta_1,\ldots,\beta_n)$ over $\cdek$},
\eeq
so that, for any fixed internal gain, maximal amplitude is achieved
simultaneously with  minimal signal duration. 
This is consistent with the notion that the most efficient cascade would 
respond with sharp (high-peaked and fast) output signals.
In the limit, this notion can be regarded as an ``instantaneous response'' 
($\sigma\approx0$) coupled with ``infinite signal amplitude'' 
($\ap\approx\infty$), which is, of course, not biologically viable. 
A realistic solution to problem (P) does exist, and is stated 
in Theorem~\ref{th-optima}.

Since the signal duration depends only on
the cascade length and the ``off'' rates, $\beta_i$,(besides the input term), we 
expect the ``on'' rates, $\alpha_i$, 
to play a small role in maximizing the efficiency of the output
response.
So, for addressing the problem (P), we will consider two different 
assumptions on the available knowledge on the $\alpha_i$: 
either (a) all the $\alpha_i$ have an equal, fixed value, $\alpha$; 
or (b) the product of the $\alpha_i$ is known, at some fixed $\alpha_P$.
We will also assume that the ``leak'' factor $\ell$ is fixed, since
this parameter was added artificially and may be adjusted independently.

\begin{figure}[h]
\epsfysize=5cm
\epsfxsize=6cm
\leftline{\hskip1cm\epsffile{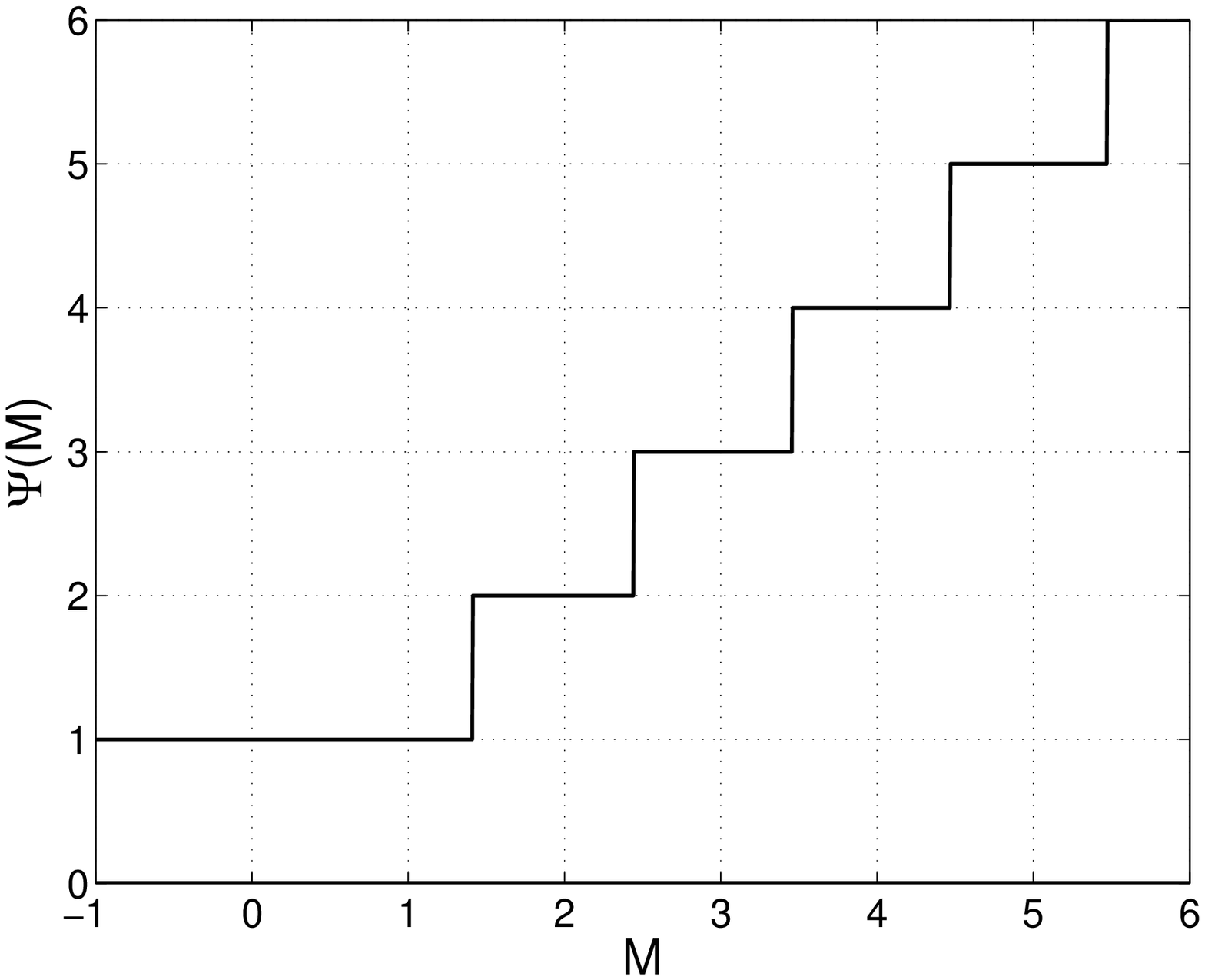}}
\vskip-5cm
\epsfysize=5cm
\epsfxsize=6cm
\rightline{\epsffile{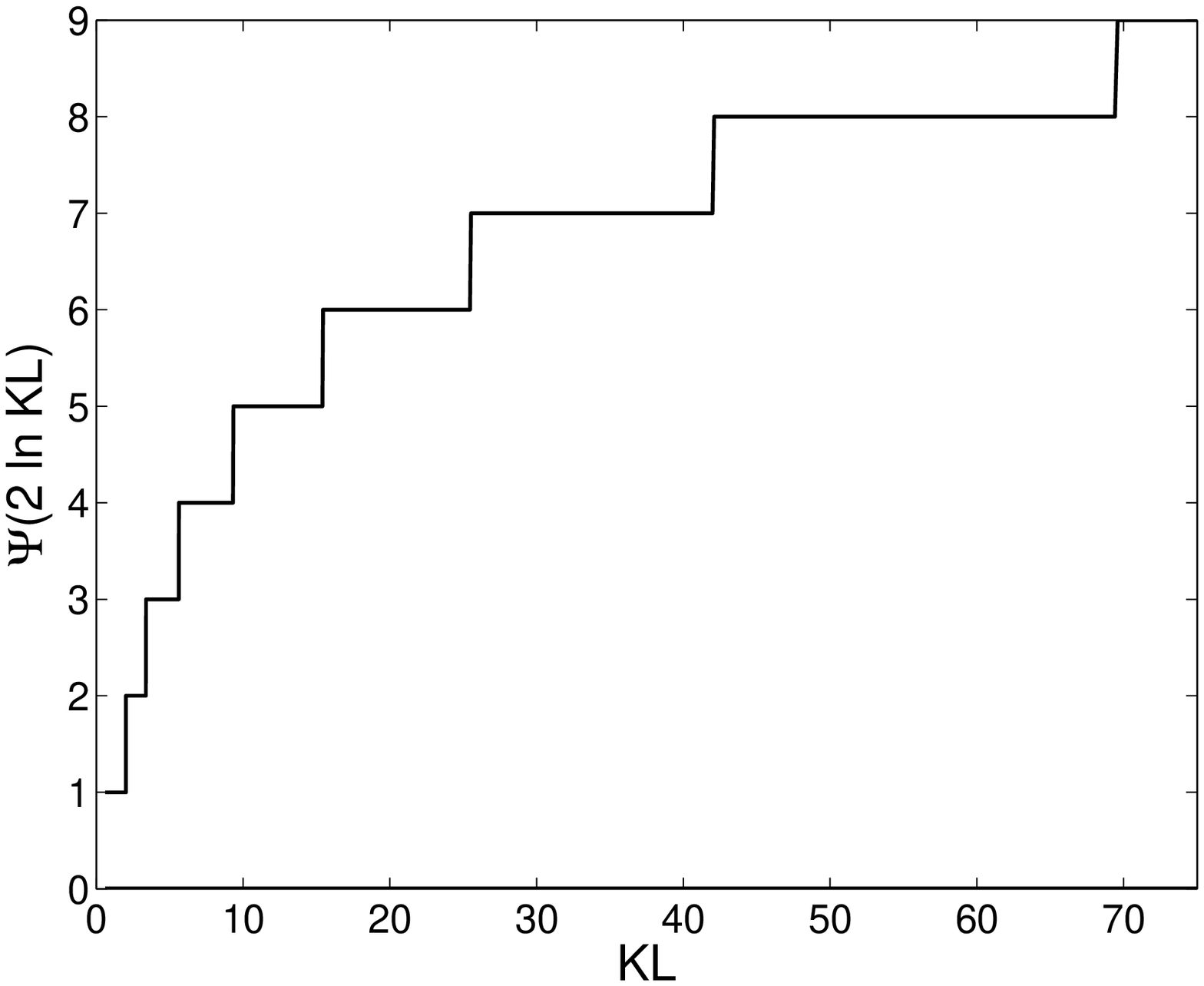}\hskip1cm}
\caption{Left: the function $\Psi(M)$.
Right: the function $\Psi(2\ln K\ell)$. Note that, for a given 
  gain $K$ and leak factor $\ell$, the optimal length is given by 
  the integer platform corresponding to the product $K\ell$. }
\label{fig-gamma}
\end{figure}

Before stating the main Theorem, we need to introduce some notation.
Define the function $f:(1,\infty)\to(0,\infty)$ to be
\beq
   f(k) = k^2\ \left[\left(1+\frac{1}{k}\right)\,
                 \ln\left(1+\frac{1}{k}\right)-\frac{1}{k}\right].
\eeq
Some properties of this function are stated in
Appendix~\ref{sec-properties}.
For any real number $M\geq1$, define
\beq
  &&  \flr{M} = \mbox{ largest integer less than or equal to $M$}, \\
  &&  \clg{M} = \mbox{ least integer greater than $M$},
\eeq
which are also known as, respectively, the ``floor'' and
``ceiling'' functions of $M$. 
Observing that any real number $M\geq1$, can be written as the sum of 
its integral and fractional parts:
\beq
   M = \flr{M} + \delta_M,
\eeq
where $\delta_M\in[0,1)$, 
define the function $\Psi:(-\infty,\infty)\to\N$,
which is plotted in Figure~\ref{fig-gamma},
\beq
   \Psi(M) = 
   \left\{\begin{array}{ll}
          1,    &    M\leq1 \\
       \flr{M}, & M>1,\mbox{ and }\delta_M\leq f(\flr{M}) \\
       \clg{M}, & M>1,\mbox{ and }\delta_M>f(\flr{M}).
    \end{array}\right.
\eeq
This is a step function where the ``jump'' between steps always occurs 
in an interval between two integers, say $k$ and $k+1$, at a point
that depends on the number $k$.
In particular, since the function $f$ is strictly increasing and
takes values in the interval $(2\ln2-1,1/2)$ (see
Appendix~\ref{sec-properties}), it follows that in some cases 
only the fractional part of the number $M$ affects the location of the
``jump'' discontinuity:
\beq
  && 0\leq \delta_M <2\ln2-1, \ \ \ \Psi(M)=\flr{M},\\
  && \frac{1}{2}< \delta_M <1,  \ \ \ \Psi(M)=\clg{M},
\eeq
while for the other cases, $2\ln2-1<\delta_M<0.5$, the choice depends
also on the integral part of $M$. 

\bt{th-optima}
Let $K>0$ and $\ell>0$ be fixed real numbers. Let $\cdek$ be the set of all 
cascades\rf{eq-in-out-sys} with internal gain $K$, as defined above. Then
\bit
\item[] 1.\ For each fixed $n=N\in\N$, the elements
$\cde=(N,\alpha_1,\ldots,\alpha_N,\beta_1,\ldots,\beta_N)\in\copt(\ell,R)$ 
satisfy $\beta_i=\beta$, for all $i=1,\ldots,N$, where 
\beq
  \beta = \left(\frac{\alpha_1\cdots\alpha_N}{K\ell}\right)^{\frac{1}{N}};
\eeq
\item[] 2(a).\ Any element $\cde\in\copt(\ell,R)$ of the form
$\cde=(n,\alpha,\ldots,\alpha,\beta_1,\ldots,\beta_n)$ 
satisfies  
\beq
    n = \Psi(2\ln K\ell)\ \ \mbox{ and }\ \ 
    \beta_i = \beta = 
    \alpha\left(\frac{1}{K\ell}\right)^{\frac{1}{n}}
\eeq
\item[] 2(b).\ Any element $\cde\in\copt(\ell,R)$ of the form
$\cde=(n,\alpha_1,\ldots,\alpha_n,\beta_1,\ldots,\beta_n)\in\copt(\ell,R)$ 
with $\alpha_1\cdots\alpha_n=\alpha_P$ satisfies
\beq
    n = \Psi\left(2\ln \frac{K\ell}{\alpha_P}\right)
    \ \ \mbox{ and }\ \ 
    \beta_i = \beta =
    \left(\frac{\alpha_P}{K\ell}\right)^{\frac{1}{n}} \ .
\eeq
\eit
\et
Before presenting the proof of the Theorem, some remarks on the
    interpretation of points 1 and 2(a), 2(b). 
The first part of the result is consistent with the observation that the ordering of
the amplification or dampening single steps within the cascade does not
influence the final output signal (also observed in~\cite{hnr2002}).

\begin{figure}[h]
\epsfysize=5cm
\epsfxsize=6cm
\leftline{\hskip1cm\epsffile{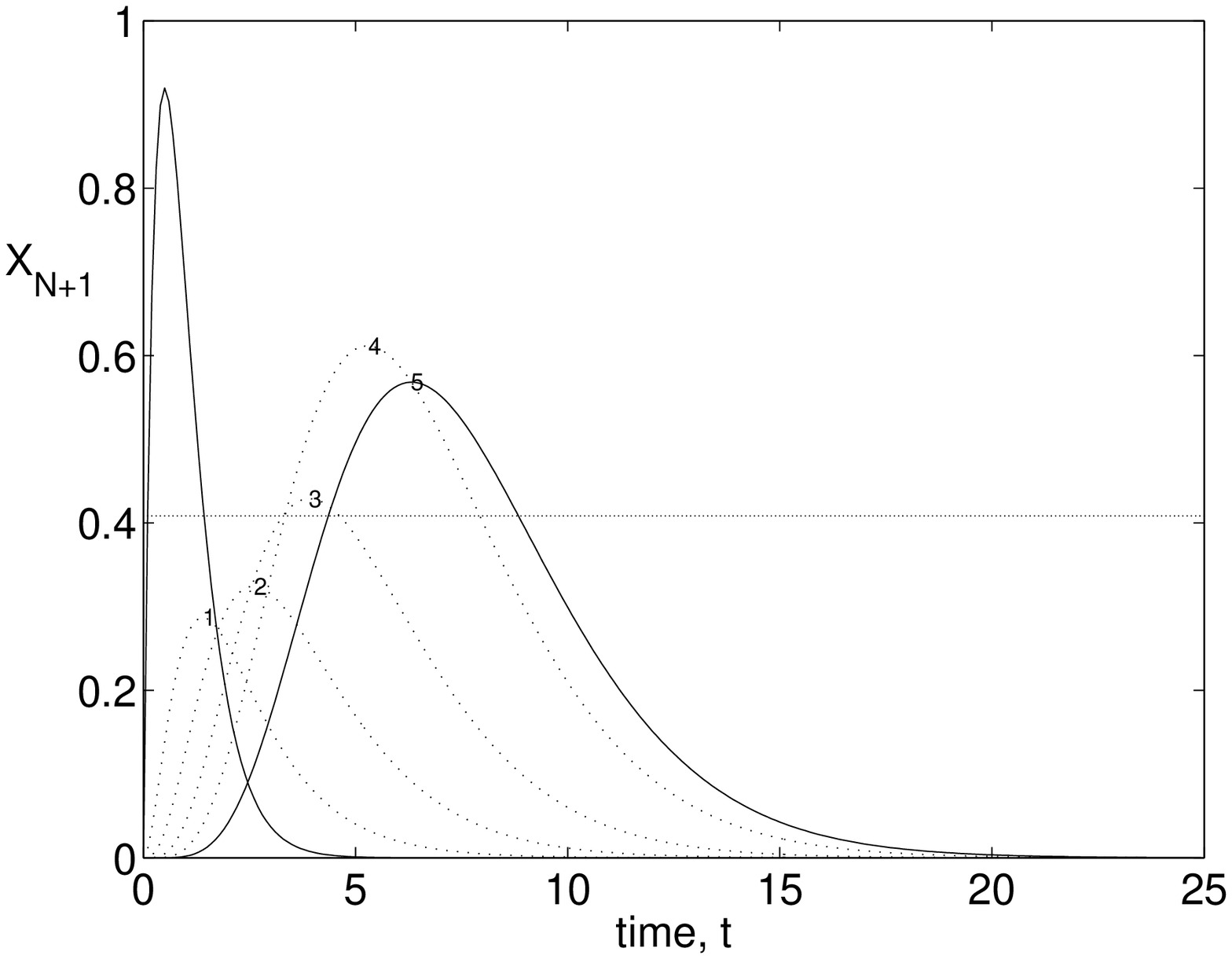}}
\vskip-5cm
\epsfysize=5cm
\epsfxsize=6cm
\rightline{\epsffile{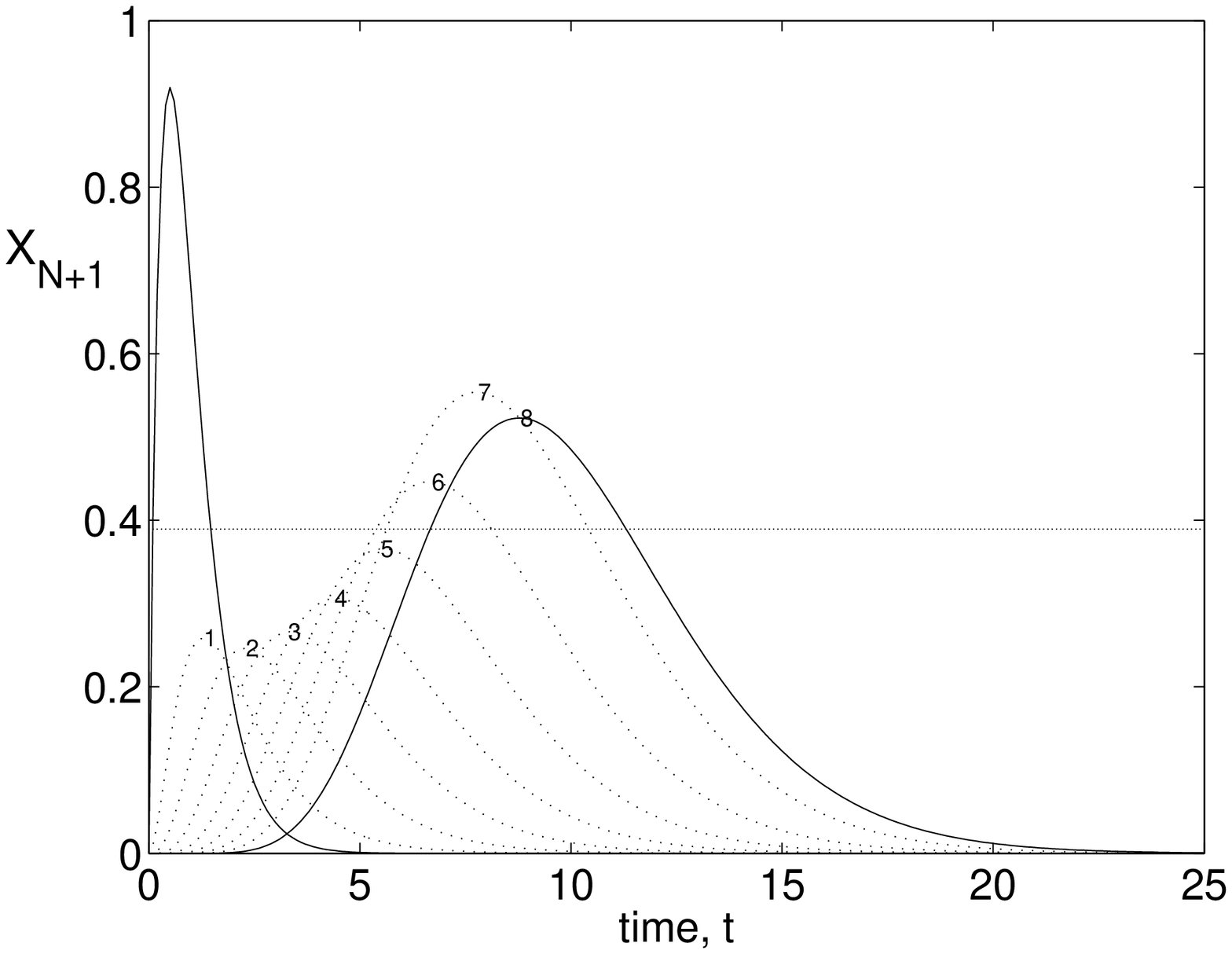}\hskip1cm}
\caption{Signal transduction cascade with $R(t)=5te^{-2t}$,
with $K=8$, $\ell=1$, $\alpha_i=1.2$. The horizontal lines represent $\ap$.
Left (optimal case): 
$n=4$, $\beta_i\approx0.714$, $i=1,\ldots,n$, $\ap=0.409$, and $\sigma_0=3.059$.
Right: $n=7$, $\beta_i\approx0.892$, $i=1,\ldots,n$, $\ap=0.389$, and $\sigma_0=3.210$.}
\label{fig-n=4-7}
\end{figure}

The second part of the Theorem shows that indefinitely increasing the
cascade's length will not increase amplification. 
In fact, there is an optimal length for the cascade that provides 
both maximum signal amplitude and duration. 
A similar observation was mentioned in~\cite{hnr2002}, and 
our Lemma~\ref{lm-Rindep} and Theorem~\ref{th-optima} characterize 
the conditions for achieving this optimization. 
For each gain $K$ and leak factor $\ell$, this optimal
length is easily read out from Figure~\ref{fig-gamma}.
For instance, a cascade with a 6 to 9-fold gain (and $\ell=1$), 
is seen to have an optimal length of 4 steps.
Figure~\ref{fig-n=4-7} illustrates Theorem~\ref{th-optima}, for an
8-fold cascade gain. The Figure shows the results of two simulations
of system\rf{eq-in-out-sys}, both with input $R(t)=5te^{-2t}$, but
different lengths of the cascade. The various curves represent $R$, the 
concentrations of each kinase $X_i$, $i=1,\ldots,n$, and the output $X_{n+1}$. 
It is clear that, for the non-optimal $n=7$,
the output's amplitude decreases and the signal duration increases. 
Note that the output curve $X_8$ is more spread out across time and its maximum
value is lower, than for the optimal case.

Theorem~\ref{th-optima} can be proved by successively solving the two
optimization problems:
\bit
\item[(P1)]  For each fixed $n$, minimize $\sigma_0$, 
over all possible choices of $\beta_1,\ldots,\beta_n\in(0,\infty)$, subject to
$\dnorm{\hat G}_\infty=K$.
\item[(P2)]  Minimize $\sigma_0$, over all possible choices of 
$n\in\N$ and $\beta_1,\ldots,\beta_n\in(0,\infty)$, subject to $\dnorm{\hat G}_\infty=K$.
\eit
Recall that we are assuming that either (a) all the $\alpha_i$ have an 
equal, fixed value, $\alpha$; 
or (b) the product of the $\alpha_i$ is known, at some fixed
$\alpha_P$.
The solution of (P1) is equal for both cases, but the solution of (P2)
is slightly different for (a) or (b).
Thus, problem (P1) is part 1 and (P2) is the part 2 of the Theorem. As we will see,
the solution of (P1) greatly simplifies the proof of (P2).

\subsection{Solving (P1): proof of part 1 of Theorem~\ref{th-optima}}
\label{ssc-solve-P1}
Given a cascade of length $n$, this problem consists of finding a set of $n$ parameters
$
   \bar\beta_1,\ldots,\bar\beta_n
$
such that the function $\sigma_0$ attains a minimum value 
at $\bar\beta_i$, $i=1,\ldots,n$, i.e.,
\beq
    \frac{1}{\bar\beta_1^2}+\frac{1}{\bar\beta_2^2}+\cdots+\frac{1}{\bar\beta_{n-1}^2}
    \leq
    \frac{1}{\beta_1^2}+\frac{1}{\beta_2^2}+\cdots+\frac{1}{\beta_{n-1}^2}
\eeq
for every $\beta_1,\ldots,\beta_n$ such that $\dnorm{\hat G}_\infty=K$:
\beq
    \dnorm{\hat G}_\infty=\frac{1}{\ell}\frac{\alpha_1\cdots\alpha_n}
                          {\beta_1\cdots\beta_n}=K
    \ \ \Leftrightarrow\ \ 
     K\ell\beta_1\cdots\beta_n - \alpha_1\cdots\alpha_n = 0.
\eeq
For simplicity, rescale the values to $B_i=1/\beta_i^2$, and observe
that
\beq
   \frac{1}{B_1\cdots B_n} = ( \beta_1\cdots\beta_n )^2 
   = \left(\frac{\alpha_1\cdots\alpha_n}{K\ell}\right)^2
\eeq
Then, the problem consists of minimizing the function:
\beq
   F(B_1,\ldots,B_{n-1})= B_1 + \cdots + B_{n-1} + \frac{Q}{B_1\cdots B_{n-1}} 
\eeq
over all possible choices of $B_i>0$, $i=1,\ldots,n-1$,
where
\beq
   Q = \left(\frac{K\ell}{\alpha_1\cdots\alpha_n}\right)^2.   
\eeq
In Appendix~\ref{sec-betas} we show that the solution to this
optimization problem is
\beq
   B_i=Q^{\frac{1}{n}},\ \ \ i=1,\ldots,n-1,
\eeq
which also implies: 
\beq
   B_n = \frac{Q}{Q^{\frac{n-1}{n}}} = Q^{\frac{1}{n}}.
\eeq
So, the choice of the ``off'' rate constants that minimizes $\sigma_0$ is 
to have $\beta_1=\beta_2=\cdots=\beta_n=\bar\beta$, with
\beq
   \bar\beta = \frac{1}{\sqrt{B_n}}= 
   \left(\frac{\alpha_1\cdots\alpha_n}{K\ell}\right)^{\frac{1}{n}},
\eeq
as we wanted to show.
$\qedb$

\subsection{Solving (P2): proof of part 2 of Theorem~\ref{th-optima}}
\label{ssc-solve-P2}
To solve the more general problem, we first show how its statement can
be simplified.
Given the value of $\alpha$ (respectively, $\alpha_P$), suppose that we have 
found a solution of (P2), i.e., an integer $n^*$  and a
set of constants $\beta^*_i$, $i=1,\ldots,n^*$ satisfying
\beqn{eq-min-sigma-nbar}
    \sigma_0(n^*,\beta^*_1,\ldots,\beta^*_{n^*})\leq 
    \sigma_0(n,\beta_1,\ldots,\beta_n)
\eeqn
for any other cascade 
$\cde=(n,\alpha_1,\ldots,\alpha_n,\beta_1,\ldots,\beta_n)$ 
with $\alpha_i=\alpha$, $i=1,\ldots,n$ (respectively, 
$\alpha_1\cdots\alpha_n=\alpha_P$).

We have already showed that
\beqn{eq-min-sigma}
    \sigma_0(n^*,\bar\beta^*,\ldots,\bar\beta^*)\leq 
    \sigma_0(n^*,\beta^*_1,\ldots,\beta^*_{n^*})
\eeqn
with 
\beq
    \bar\beta^*=\left( \frac{\alpha_1\cdots\alpha_{n^*}}{K\ell}
                \right)^{\frac{1}{n^*}}
\eeq
and we know this choice yields the unique minimum of $\sigma_0$ 
for a fixed length $n$. So, it follows that the solution of (P2) must also satisfy
\beq
    \beta^*_i=\bar\beta^*,\ \ \ i=1,\ldots,n^*.
\eeq
This observation allows us to simplify the statement of problem (P2),
and look only for solutions where all $\beta_i$'s are equal:
\bit
\item[(P2)']  Minimize $\sigma_0(n,\beta,\ldots,\beta)=n/\beta^2 $, 
over $n\in\N$ and $\beta\in(0,\infty)$,
subject to $(\alpha/\beta)^n=K\ell$.
\eit
From the constraint $\dnorm{\hat G}_\infty=K$ we have 
\beq
\mbox{case 2(a): } &&
   \left(\frac{\alpha}{\beta}\right)^n=K\ell  \ \ \Leftrightarrow\ \ 
   \beta = \alpha\ \left(\frac{1}{K\ell}\right)^{\frac{1}{n}}
   \ \ \Rightarrow\ \ 
   \sigma_0(n,\beta(n)) = \frac{1}{\alpha^2}\ n\ (K\ell)^{\frac{2}{n}}.
\\ &&  \\
\mbox{case 2(b): }&&
   \frac{\alpha_P}{\beta^n}=K\ell  \ \ \Leftrightarrow\ \ 
   \beta = \left(\frac{\alpha_P}{K\ell}\right)^{\frac{1}{n}}
   \ \ \Rightarrow\ \
   \sigma_0(n,\beta(n)) = n\ \left(\frac{K\ell}{\alpha_P}\right)^{\frac{2}{n}}.
\eeq
In either case, to solve the problem, it is enough to minimize the function
$\ln[\sigma_0(n,\beta(n))]$:
\beq
   F(n,M) = \ln n +\frac{1}{n}\; M 
\eeq
over $n\in\N$, where $M$ is a positive constant with value either
\beqn{eq-Mvalue}
   M &=& 2 \ln\; K\ell,\ \ \mbox{ for case 2(a) } \\ 
     & &\nonumber\\
   M &=& 2 \ln\; \frac{K\ell}{\alpha_P},\ \ \mbox{ for case 2(b) }.
\eeqn
For a fixed $M$, let the minimizer of $F(n,M)$ over $n\in\N$ be 
\beq
   n^*(M):= \{ n\in\N:\ F(n,M)\leq F(n',M),\ \mbox{ for every }\ n'\in\N \}, 
\eeq
which is given by Lemma~\ref{lm-minFn} (Appendix~\ref{sec-nstar}):
\beq
    n^*(M)=\Psi(M).
\eeq
Thus, for part 2(a) of the Theorem we have 
$n=n^*(2 \ln\; K\ell)=\Psi(2 \ln\; K\ell)$, and for part 2(b) we
have $n=\Psi(2 \ln\; K\ell/\alpha_P)$.
The value $\beta$ is given according to part 1.
$\qedb$

\begin{figure}[h,t]
\epsfysize=5cm
\epsfxsize=6cm
\leftline{\hskip5mm\epsffile{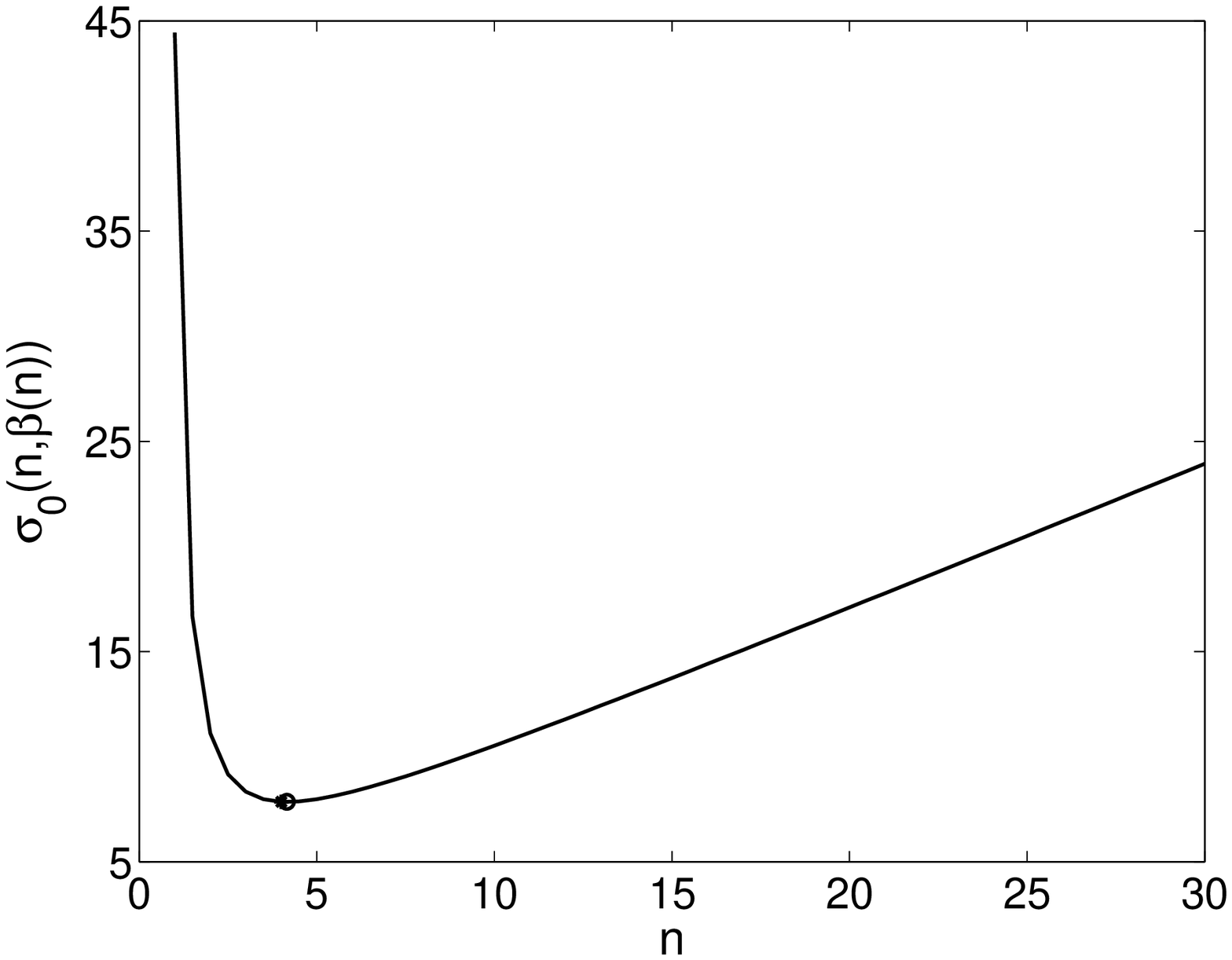}}
\vskip-5cm
\epsfysize=5cm
\epsfxsize=6cm
\hskip-5mm
\rightline{\epsffile{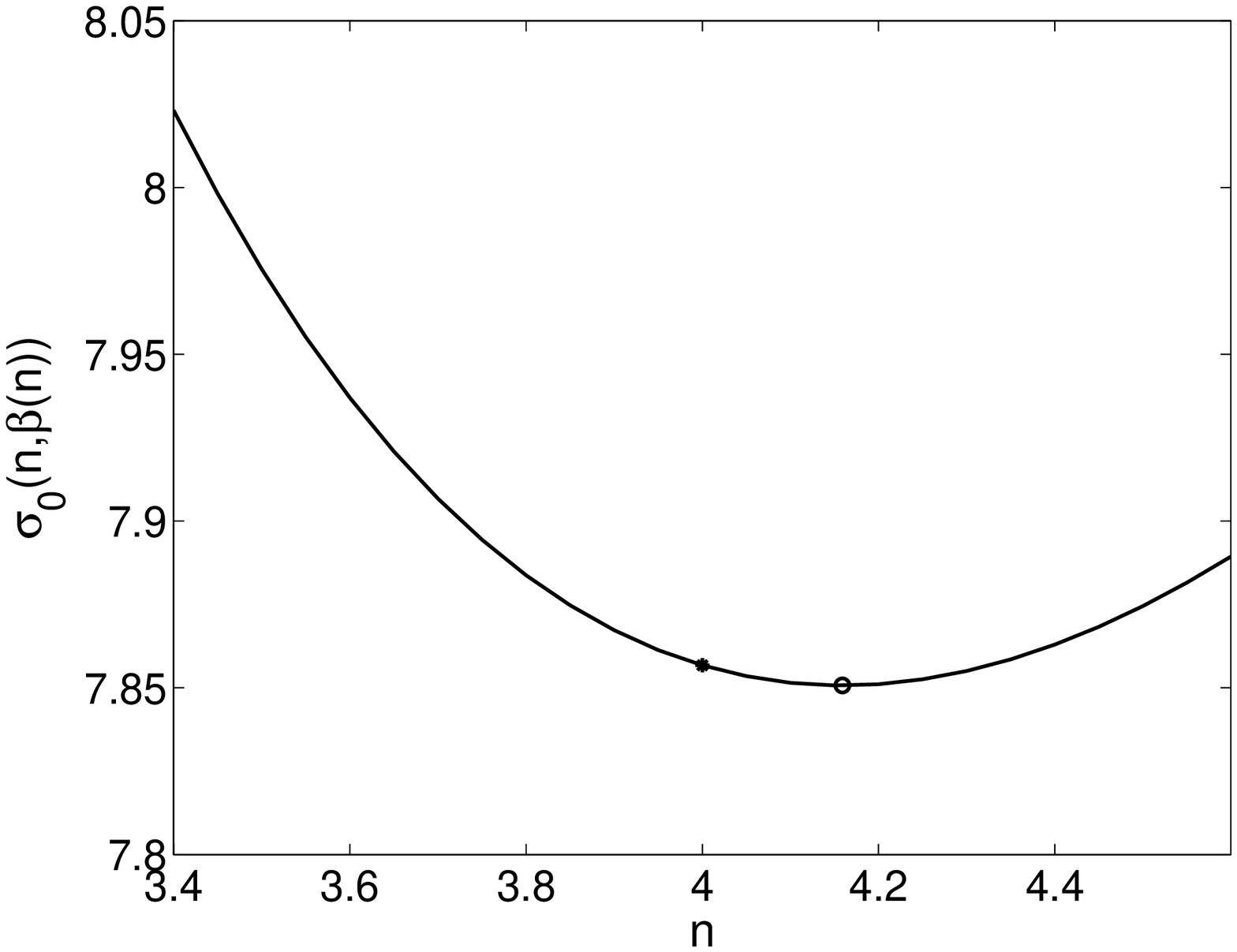}\hskip5mm}
\caption{The function $\sigma_0(n,\beta(n))$, for $K\ell=8$ and
  $\alpha=1.2$; and the points
  $(M,\sigma_0(M,\beta(M)))$ (circle) and
  $(n^*,\sigma_0(n^*,\beta(n^*)))$ (star).}
\label{fig-sigma0}
\end{figure}

\vspace{1cm}
As shown in the example of Figure~\ref{fig-sigma0}, 
evaluation of $\sigma_0$ at $n^*(M)$
yields a value which is actually quite
close to the ``true'' $\sigma_0(M,\beta(M))$.

\section{Cascades with positive feedback}
\label{sec-feedback}
In this Section we investigate the behavior of cascades under 
positive feedback. Assume that the last kinase, $X_n$, also
contributes to the activation of the first kinase: 
then the differential equation for $X_1$ includes one more term and
becomes 
\beq
   \frac{dX_1}{dt} = \alpha_1R(t)+\eps X_n -\beta_1X_1.
\eeq
We will assume that $\eps$ is small enough: 
\beq
   \beta_1\cdots\beta_n > \eps\alpha_2\cdots\alpha_n.
\eeq 
This guarantees that the cascade is stable with respect to small 
perturbations (that is, all the eigenvalues of the system's matrix 
$A$ have negative real parts, see Section~\ref{sec-stability}).

We can compute the transfer function for the system with feedback 
($\eps>0$),
just as we did in Section~\ref{sec-tf}, for a given cascade
$\cde=(n,\alpha_1,\ldots,\alpha_n,\beta_1,\ldots,\beta_n)$
and any input $R$, and leak factor $\ell>0$. We obtain:
\beqn{eq-Ghat-fdb}
   \hat G(s)=\frac{1}{s+\ell}\ \frac{\alpha_1\cdots\alpha_n}
                {(s+\beta_1)\cdots(s+\beta_n)-\eps\alpha_2\cdots\alpha_n}.
\eeqn    
The infinity norm is again obtained for the case $\omega=0$ 
($s=\jmath\omega$) (see below, at the end of this Section):
\beq
   \dnorm{\hat G}_{\infty} = \frac{1}{\ell}\ 
             \frac{\alpha_1\cdots\alpha_n}
                  {\beta_1\cdots\beta_n-\eps\alpha_2\cdots\alpha_n}.
\eeq
Computing the signaling time ($\tau$), and the signal duration
($\sigma$) and amplitude ($\ap$), we have
\beqn{eqc-tau-fdb}
    \tau_{\fbk}(\cde,\ell,R)=\frac{1}{\ell}
         +\frac{\beta_1\cdots\beta_n\sum_{i=1}^{n}\frac{1}{\beta_i}}
               {\beta_1\cdots\beta_n-\eps\alpha_2\cdots\alpha_n}
         + \left.\frac{d\ \ln\hat R}{ds}\right\rfloor_{s=0}
\eeqn
\beqn{eqc-sigma-fdb}
    \sigma_{\fbk}(\cde,\ell,R)=
    \sqrt{\frac{1}{\ell^2}
    +\frac{(\beta_1\cdots\beta_n)^2\left[
            \sum_{i=1}^{n}\frac{1}{\beta_i^2}
           +\eps\alpha_2\cdots\alpha_n
             \sum_{i\neq j}^{}\frac{1}{\beta_i\beta_j}\right] }   
               {(\beta_1\cdots\beta_n-\eps\alpha_2\cdots\alpha_n)^2}
    + q(R)},
\eeqn
\beqn{eqc-amplitude-fdb}
   \ap_{\fbk}(\cde,\ell,R)
       =\frac{1}{\sigma(\cde,\ell,R)}\
       \frac{\alpha_1\cdots\alpha_n}
          {\beta_1\cdots\beta_n-\eps\alpha_2\cdots\alpha_n}\ \dnorm{R}_2,
\eeqn
Comparison of these quantities for the models with and without
feedback leads to the following conclusions:
\bit
\item[1.] the system with feedback exhibits higher internal gain;
\item[2.] the system with feedback exhibits larger signaling time 
and signal duration $\tau_{\fbk}>\tau$ and $\sigma_{\fbk}>\sigma$.
\eit

So, for an arbitrary cascade, the existence of a positive feedback
leads to a less sharp output signal: the signal transduction
down the cascade takes a longer time, and the output signal has 
greater duration. 

On the other hand, the existence of feedback may be used to great
advantage in the design of an optimal cascade: positive feedback
(at a constant rate $\eps$) allows the cascade to be of shorter length
and still have the same maximal amplitude/minimal duration.
The results in Theorem~\ref{th-optima} are valid just as before, 
with suitable adjustements to some of the constants. Thus, we now have
\beq
   \dnorm{\hat G}_{\infty}=K \ \ \Leftrightarrow\ \ 
   \beta_1\cdots\beta_n = 
   \frac{(\alpha_1+\eps K\ell)\alpha_2\cdots\alpha_n}{K\ell},
\eeq 
and now, similarly to the proof in Section~\ref{ssc-solve-P1}, we set
\beq
   Q_{\fbk}=\left(\frac{K\ell}
             {(\alpha_1+\eps K\ell)\alpha_2\cdots\alpha_n}\right)^2,
\eeq
which leads to the optimal value for $\beta_i=\beta_{\fbk}$, $i=1,\ldots,n$
\beq
   \beta_{\fbk}=\left(\frac{(\alpha_1+\eps K\ell)\alpha_2\cdots\alpha_n}
                    {K\ell}\right)^{\frac{1}{n}}.
\eeq
To find the optimal length of the cascade with feedback, note that
\beq
  \sigma_0(n,\beta(n)) = n M_{\fbk}^{\frac{1}{n}},\ \ \mbox{ with }
  M_{\fbk}=2\ln\ \left(\frac{K\ell}
                            {(\alpha_1+\eps K\ell)\alpha_2\cdots\alpha_n}
                    \right)^{\frac{1}{n}}.
\eeq
Since $M_{\fbk}\leq M$, then also $n^*_{\fbk}(M_{\fbk})\leq n^*(M)$. 
Therefore, we conclude that, for the cascade with feedback,
\bit
\item[3.] for each fixed $n$, the value of the off rates that maximixes
$\ap$ (minimizes $\sigma_0$) over $\copt(\ell,R)$ is larger,
$
    \beta_{\fbk}>\beta;
$
\item[4.] the length of the cascade that maximixes
$\ap$ (minimizes $\sigma_0$) over $\copt(\ell,R)$ is smaller,
$
    n^*_{\fbk} < n^*.
$
\eit

These results agree with what would be expected from a
signaling pathway: indeed, the
existence of positive feedback enhances the activation at each step, 
so a larger amount of the phosphorylated kinase will be produced; to
keep this amount at a ``weak'' level, the phosphatases should increase
their activity.
On the other hand, since the amount of phosphorylated kinases
increased, a smaller number of steps is required to produce the same
signal amplitude as in the cascade with no feedback.

\vspace{0.5cm}
To compute the infinity norm $\dnorm{\hat G}_{\infty}$, we first note that 
the denominator of $\hat G(\jmath\omega)$, which we will denote by
$\texttt{den}(\hat G(\jmath\omega))$,  satisfies (by the triangle inequality):
\beq
   \abs{\texttt{den}(\hat G(\jmath\omega))}  \geq
   \abs{\jmath\omega+\ell}\ [\abs{\jmath\omega+\beta_n}\cdots\abs{\jmath\omega+\beta_n}
                           -\eps\alpha_2\cdots\alpha_n].
\eeq
Also
\beq
  \abs{\jmath\omega+\beta_n}\cdots\abs{\jmath\omega+\beta_n} 
  = \sqrt{(\omega^2+\beta_1^2)\cdots(\omega^2+\beta_n^2)}
  \geq \beta_1\cdots\beta_n,
\eeq
for every $\omega\in\R$, where the equality holds if and only if $\omega=0$.
Thus
\beq
  \abs{\texttt{den}(\hat G(\jmath\omega))} 
  \geq \ell\;[\beta_1\cdots\beta_n -\eps\alpha_2\cdots\alpha_n]
  = \texttt{den}(\hat G(0)) >0,
\eeq
where the last inequality follows from
the assumption $\beta_1\cdots\beta_n >\eps\alpha_2\cdots\alpha_n$.
Therefore, if the expression $\abs{\texttt{den}(\hat G(\jmath\omega))}$
is minimized at $\omega=0$,
then the function $\abs{\hat G(\jmath\omega)}$ is maximized at
$\omega=0$, as we wanted to show.

\section{Signal delay}
\label{sec-delay}
The frequency domain approach for linear systems also provides an
answer to certain problems involving delays and stability within a signaling cascade.
For instance, if there is delay in transmitting the signal at any step 
along the cascade, then both the amplitude and the signal duration are 
not affected. Suppose that, at each step, there is a delay $\delta_i$ in the
transmission of the signal. The differential equation becomes
\beq
   \frac{d}{d t}\pmatrix{X_1(t)\cr X_2(t)\cr \vdots\cr X_{n+1}(t)}=
   A\pmatrix{X_1(t-\delta_1)\cr X_2(t-\delta_2)\cr 
               \vdots\cr X_{n+1}(t-\delta_{n+1})}\ +\ B\;R.
\eeq
The Laplace transform of $X_i(t-\delta_i)$ is, from the properties
listed in the Appendix,
\beq
   e^{-s\delta_i}\hat X_i
\eeq
so that
\beq
     \hat X_{i+1}=\frac{\alpha_{i+1}}{s+\beta_{i+1}}
                  e^{-s\delta_i}\hat X_i.
\eeq
The transfer function becomes:
\beq
    \hat G(s)=\frac{1}{s+\ell}\frac{\alpha_1\cdots\alpha_n}{(s+\beta_1)\cdots(s+\beta_n)}
         \ e^{-s\delta_1}\cdots e^{-s\delta_{n+1}}.
\eeq
But, for an imaginary number $\jmath\omega$, 
$\abs{e^{-\jmath\omega\delta_i}}=1$, so the norm 
$\dnorm{\hat G}_{\infty}$ is unchanged; and since $e^{-s\delta_i}=1$ when
evaluated at $s=0$, the signal duration and amplitude are also unchanged. 
This is not surprising, because in a linear system, delay
simply causes a temporal translation of the signal, by a fixed amount,
without affecting amplitudes.

\section{Stability of cascades}
\label{sec-stability}
A signaling pathway is considered {\it stable} (see~\cite{hnr2002}) if
small and random perturbations (those that do not consist of biologically 
relevant inputs) are not amplified. 
So, in the presence of small perturbations, the amount of phosphorylated 
kinases should not be allowed to grow very large, and should return to the
stable state, with $X_i\approx0$, for all $i=1,\ldots,n$.
Thus, the behavior of a signaling pathway in the absence of a relevant input always 
satisfies expression\rf{eq-weak-act}, that is, $X_i\ll X_{\tot,i}$
for each $i=1,\ldots,n$, and hence its stability may be established
by analysis of the model\rf{eq-in-out-sys}.

In the absence of an input ($R(t)\equiv0$), the point
$(X_1,X_2,\ldots,X_{n+1})=(0,0,\ldots,0):={\bf 0}$ is an equilibrium point of 
system\rf{eq-in-out-sys}, and the stability of this equilibrium
determines the stability of the pathway.
The equilibrium point ${\bf 0}$ is stable if all the eigenvalues of the 
matrix $A$ have negative real parts. This is indeed the case for the
system described by equations\rf{eq-in-out-sys}.
We know that, after a perturbation, the system will always return to ${\bf 0}$.
Moreover, we can estimate that a small perturbation will also generate a small
response, since:
\beq
    \dnorm{Y_{pert}}_2 \leq \kappa\dnorm{R_{pert}}_2,
\eeq
where $\kappa$ is a constant, equal to $\dnorm{\hat G}_{\infty}$.

For signaling cascades which exhibit a lower degree of kinase specificity,
the problem of stability of the cascade (see~\cite{hnr2002}) becomes 
significant. If a kinase $X_i$ affects both the downstream kinases and some 
upstream kinase, then the eigenvalues of $A$ change, and stability is not guaranteed.
Allowing for kinase non-specificity, a resulting matrix $A$ could be of the form:
\beq
    A_\eps=\pmatrix{-\beta_1 & \eps_{12} & 0 & \cdots & 0 & \eps_{1n} & 0 \cr
               \alpha_2 & -\beta_2 & 0 & \cdots & 0 & 0 & 0\cr
               \eps & \alpha_3 & -\beta_3 & \cdots & 0 & 0 & 0\cr
               \eps & \eps & \alpha_4 & \cdots & 0 & 0 & 0\cr
               \vdots &  & & \ddots &  & \vdots \cr
               \eps & \eps & \eps & \cdots &\alpha_n & -\beta_n & 0\cr
                0 & 0 & 0 & \cdots & 0 & 1 & \ell\cr
               },
\eeq
and the signaling pathway is stable if
\beq
   \mbox{all the eigenvalues of $A_\eps$ have negative real parts.}
\eeq
Some relevant easy-to-compute examples are:
\bit
\item suppose that each kinase $i$ is only allowed to activate its 
downstream kinases ($\eps_{1n}=0$, $\eps_{12}=0$ and $\eps\neq0$); 
in this case 
it is not surprising that stability is not affected at all, because this situation
corresponds to a lower triangular matrix, again with eigenvalues $-\beta_i$;
\item suppose that there exists feedback from the last activated kinase
to the first kinase ($\eps_{1n}=\eps_0$, $\eps_{12}=0$ and $\eps=0$);
in this case, if $\eps_0$ satisfies
\beq
   \beta_1\cdots\beta_n>\eps_0\alpha_2\cdots\alpha_n
   \ \ \Leftrightarrow\ \ 
   \eps_0<\frac{\alpha_2\cdots\alpha_n}{\beta_1\cdots\beta_n},
\eeq
then the eigenvalues of the matrix $A_{\eps}$ all have negative real
parts and the new cascade is stable.
To prove this, suppose that there exists an eigenvalue of $A_{\eps}$
with positive real part, that is, a complex number $\lambda$ such that
\beqn{eq-lambda}
   \lambda=\lambda_{re}+\jmath\lambda_{im}, \ \ \ \ 
   \lambda_{re}\geq0,
\eeqn
and
\beqn{eq-ldet}
    \mbox{det}(A_\eps-\lambda I)=
    (-1)^n\ell[(\lambda+\beta_1)\cdots(\lambda+\beta_n) -\eps_0\alpha_2\cdots\alpha_n] =0.
\eeqn
Then\rf{eq-lambda} implies  $\abs{\lambda+\beta_i}\geq\beta_i$, for
$i=1,\ldots,n$ and so
\beq
   \abs{(\lambda+\beta_1)\cdots(\lambda+\beta_n)
     -\eps_0\alpha_2\cdots\alpha_n}
   &\geq&
   \abs{(\lambda+\beta_1)\cdots(\lambda+\beta_n)}-\eps_0\alpha_2\cdots\alpha_n \\
   &\geq& \beta_1\cdots\beta_2-\eps_0\alpha_2\cdots\alpha_n>0,
\eeq
which contradicts equation\rf{eq-ldet}.
\eit

\section{Conclusions}
By modeling weakly activated signal transduction cascades  as
linear systems and applying techniques from control systems theory one can
identify the cascade's input-to-output transfer function and
internal gain. Based on these properties, the concepts of signal
duration, signaling time and signal amplitude may be defined in 
an intuitive and general form, for any input signal.

Our analysis shows that, for linear cascades, signal amplitude and 
duration are, respectively, maximized and minimized simultaneously.
So, a cascade can respond with signals that are both fast and exhibit
high amplification.
To achieve the highest amplification and the shortest duration response, 
the cascade should have all off rates equal to some value $\beta$.

We also show that, for each fixed internal gain, there are {\it finite 
values} for the length of the cascade and the off constants that 
simultaneously maximize (resp., minimize) the signal amplitude 
(resp., signal duration). 
To achieve these optimal conditions, the optimal length should be
given by the well defined step function $\Psi$.
This function $\Psi$ depends only on, and increases logarithmically
with, the internal gain of the system. 
The off constants should all have the same value $\beta$. 
This optimal value $\beta$ 
depends on the internal gain and the length of the system.
 
In addition, our analysis shows that a positive feedback term 
on the cascade enhances the optimal design, by allowing the same
signal amplitude and duration to be achieved with a shorter length
and higher off rates.

Finally, other issues, such as delay at each phosphorylation step, 
and the stability of the signaling pathway when there is a high degree of
non-specificity among the kinases, are also naturally examined within
this framework. The stability of the zero steady-state of the cascade
with respect to small perturbations is established by checking that
the eigenvalues of the matrix $A$ all have negative real parts.

\appendix
\section{Properties of function $f(k)$}
\label{sec-properties}

The function  $f:(1,\infty)\to(0,\infty)$ 
\beq
   f(k) = k^2\ \left[\left(1+\frac{1}{k}\right)\,
                 \ln\left(1+\frac{1}{k}\right)-\frac{1}{k}\right]
\eeq
has the following properties:
\bit
\item[1.] $f$ is strictly increasing;
\item[2.] $f(1)=2\ln2-1\approx0.386$ and $\lim_{k\to\infty}f(k)=1/2$.
\eit
To prove property 1, notice that another expression for $f$ is 
$f(k)=k[(k+1)\ln\ (k+1)/k \ -\ 1]$, and compute the first and second derivatives:
\beq
  \frac{d\ f}{d k} &=& (2k+1)\ \ln\;\frac{k+1}{k} - 2 \\
  \frac{d^2\ f}{d k^2} &=& (2k+1)\  \ln\;\frac{k+1}{k}+\frac{2k+1}{k(k+1)}.
\eeq
It is clear that the second derivative is always positive, and hence
the first derivative is strictly increasing. Since
$df/dk(1)=3\ln2-1>0$, it follows that the first derivative
is also always positive and therefore the function $f$ is strictly
increasing.

To prove property 2, the value $f(1)$ is straightfoward, and for the
limit as $k\to\infty$, it is easier to consider $x=1/k$ and compute:
\beq
  \lim_{k\to\infty}f(k)=\lim_{x\to0}f(1/x)
 =\lim_{x\to0}\frac{(1+x)\ln(1+x) - x}{x^2}= \frac{0}{0}.
\eeq
This indeterminacy can be solved by twice applying L'H\^opital's rule: 
\beq
   \mbox{...first time:}\ \ \ \  &&
   \frac{\ln(1+x) +(1+x)\frac{1}{1+x} - 1}{2x} \to\frac{0}{0},\ \ \
                                              \mbox{as $x\to0$} \\ 
   && \\
   \mbox{...second time:}\ \ \ \ &&
   \frac{\frac{1}{1+x}}{2} \to \frac{1}{2},\ \ \
                               \mbox{as $x\to0$}.
\eeq

\begin{figure}[h,t]
\epsfysize=5cm
\epsfxsize=6cm
\centerline{\epsffile{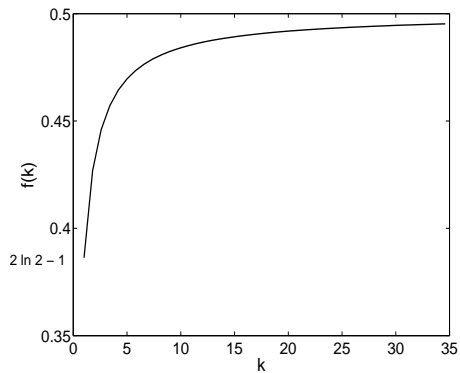}}
\caption{The function $f(k)$.}
\label{fig-f(k)}
\end{figure}

\section{Minimization of $\sigma_0$}
\label{sec-betas}

Let $Q$ be a positive real number and $n>2$ an integer.
Consider the function $F:(0,\infty)^{n-1}\to(0,\infty)$ given by
\beq
   F(B_1,\ldots,B_{n-1})= B_1 + \cdots + B_{n-1} + \frac{Q}{B_1\cdots B_{n-1}}.
\eeq

\bl{lm-betas}
The choice of $B_i>0$, $i=1,\ldots,n-1$ that
minimizes the function $F$ is: $B_i=Q^{1/n}$, $i=1,\ldots,n-1$.
\el

\bpr
First, we claim that the search for a point $(B_1,\ldots,B_{n-1})$ 
where $F$ is minimized  can be constrained to the compact set:
\beqn{eq-bdyset}
    \left[ \frac{1}{n^{n-1}}\;Q^{1/n},\ n\;Q^{1/n} \right]^{n-1}.
\eeqn
To justify the upper bound of the interval, observe that
\beqn{eq-min-point}
    F(Q^{1/n},\ldots,Q^{1/n}) = n\;Q^{1/n} 
\eeqn
and that, for any $j=1,\ldots,n-1$,  
\beqn{eq-bdy1}
    B_j\geq n\;Q^{1/n} \ \ \Rightarrow\ \ F(B_1,\ldots,B_{n-1})>n\;Q^{1/n}.
\eeqn
So, it is enough to look for a minimum of $F$ in the region
$B_i< n\;Q^{1/n}$, $i=1,\ldots,n-1$
(because inside this region there is at least one point --
equation\rf{eq-min-point} -- where $F$ has a lower value than  
anywhere outside of this region).

To justify the  lower bound, suppose that
$B_j\leq\frac{1}{n^{n-1}}Q^{1/n}$, for some $j=1,\ldots,n-1$. 
Then using the already established upper bounds
\beqn{eq-bdy2}
   F(B_1,\ldots,B_{n-1})
     >\frac{Q}{B_1\cdots B_j\cdots B_{n-1}} 
 \geq \frac{Q}{\frac{1}{n^{n-1}}Q^{1/n}\ [n\;Q^{1/n}]^{n-2}}
     = n\; Q^{1/n},
\eeqn
and similarly we conclude that it is enough to look for a minimum 
of $F$ in the region $B_i>\frac{1}{n^{n-1}}Q^{1/n}$, for $i=1,\ldots,n-1$.

The function $F$ is continuous, in fact differentiable, in the compact 
set\rf{eq-bdyset}, and so $F$ has (absolute) maximum and minimum values in this
set. The maximum and minimum may be attained either at a critical
point of $F$, or at the boundary points of\rf{eq-bdyset}. Equations\rf{eq-bdy1}
and\rf{eq-bdy2} show that the minimum is not attained at any of the
boundary points. So the minimum will be attained at an interior
point of the set\rf{eq-bdyset}, which must also be a critical point of
$F$. The critical points of $F$ are given by:
\beq
   \frac{dF}{dB_j} = 0 \ \ \Leftrightarrow\ \ 
   1 - \frac{1}{B_j}\frac{Q}{B_1\cdots B_{n-1}},\ \ j=1,\ldots,n-1,
\eeq
or, equivalently,
\beq
   B_j = \frac{Q}{B_1\cdots B_{n-1}} = B_*,\ \ j=1,\ldots,n-1,
\eeq
where  $B_*$ satisfies
\beq
   1 - \frac{1}{B_*}\frac{Q}{B_*^{n-1}} =0
   \ \ \Leftrightarrow\ \ 
   B_* = Q^{1/n}.
\eeq
Thus, there exists a unique critical point of $F$, $(B_*,\ldots,B_*)$,
which indeed belongs to the compact set\rf{eq-bdyset}. 
By the discussion above, this point must be the minimizer of $F$, as
we wanted to show.
\epr

\section{Minimization of $F(n,M)$}
\label{sec-nstar}

For a fixed $M$, let the minimizer of $F(n,M)$ over $n\in\N$ be 
\beq
   n^*(M):= \{ n\in\N:\ F(n,M)\leq F(n',M),\ \mbox{ for every }\ n'\in\N \}, 
\eeq

\bl{lm-minFn}
Let $M$ be any fixed real number. Then $n^*(M)=\Psi(M)$.  
\el

\bpr
Since $n^*(M)$ is the minimizer of $F(n,M)$ over the (positive)
natural numbers, we start by computing the derivative of $F(n,M)$:
\beq
   \frac{d\, F}{dn}(n,M) = \frac{1}{n} -\frac{1}{n^2} M
                         = \frac{1}{n^2}\ [n - M].
\eeq
We consider two distinct cases:
\bit
\item Case $M\leq0$
\beq
   \frac{d\, F}{dn}(n,M) > 0,\ \ \mbox{ for all $n\geq1$},
\eeq
so $F(\cdot,M)$ is a strictly increasing function and thus its minimizer
over $\N$ is the smallest natural number, i.e., $n^*(M)=1$.
\item Case $M>0$
\beq
   \frac{d\, F}{dn}(n,M) =0\ \ \Leftrightarrow\ \ n=M, 
\eeq
and the derivative is negative for $n<M$ and positive for $n>M$: in
other words, the function $F$ has indeed a minimum at $n=M$.
However, in general, $M$ is not an integer, so it cannot be a solution
to our minimization problem. We should choose 
\beq
   n^*(M) = \left\{\begin{array}{ll}
                  1,    & M\leq1\\
               \flr{M}, & M>1, \mbox{ and } F(\clg{M},M)\geq F(\flr{M},M)\\
               \clg{M}, & M>1, \mbox{ and } F(\clg{M},M)< F(\flr{M},M). 
            \end{array}\right.
\eeq
Note that we pick $n^*=1$ whenever $M\leq1$, since a cascade of length
zero is meaningless. 
\eit
To further analyze this condition, observe that we can write, for $M>1$,
\beq
    M=k+\delta,\ \flr{M}=k,\ \ \clg{M}=k+1
\eeq 
where $k\geq1$ is the integral part of $M$ 
and $\delta\in[0,1)$ is the fractional part of $M$. 
Now, the point $\delta$ for which $n^*$ ``jumps'' from $\flr{M}$ to $\clg{M}$
can be found by setting 
\beq
   0 &=& F(\clg{M},M) - F(\flr{M},M) = F(k+1,k+\delta) - F(k,k+\delta) \\
     &=&  \ln (k+1) +\frac{1}{k+1} (k+\delta) - \ln k - \frac{1}{k}(k+\delta).
\eeq
Simplifying this equation we obtain:
\beq
     \ln\frac{k+1}{k} - \frac{k+\delta}{k(k+1)} = 0
     \ \ &\Leftrightarrow&\ \ 
     \delta = k(k+1)\,\ln\frac{k+1}{k} -k  \\ 
     \ \ &\Leftrightarrow&\ \ 
     \delta = k^2\left[\frac{k+1}{k}\,\ln\frac{k+1}{k} -\frac{1}{k}\right]\\
      \ \ &\Leftrightarrow&\ \ 
     \delta = k^2\left[\left(1+\frac{1}{k}\right)\,
                 \ln\left(1+\frac{1}{k}\right)-\frac{1}{k} \right]
            = f(k).
\eeq
Analysis of this function (see Appendix~\ref{sec-properties}), shows
that $f$ is positive and strictly increasing, so we have
\beq
   F(\clg{M},M) - F(\flr{M},M) \geq 0 \ \ \Leftrightarrow\ \ 
   f(\clg{M}) - \delta \geq 0.
\eeq
Therefore, we should choose 
\beq
  n^*(M) = \left\{\begin{array}{ll}
                  1,    & M\leq1\\
              \flr{M}, & M>1, \mbox{ and }  \delta \leq f(\flr{M})\\
              \clg{M}, & M>1, \mbox{ and } \delta >f(\flr{M}). 
           \end{array}\right.
\eeq
This proves the Lemma.
\epr

\section{Dictionary: Laplace transforms and transfer functions}
\label{sec-fd}

For further details about these topics see, for instance,~\cite{rota}
and~\cite{control},~\cite{control2},~\cite{mct}.

\subsubsection*{Laplace transforms}
For a function $X:(-\infty,\infty)\to\R^n$ 
(with $\abs{X(t)}\leq c e^{kt}$, for all $t$, for some positive constants $c$, $k$),  
the Laplace transform is another function $\hat X:\mathcal R\to\C^n$  defined as 
\beq
    \hat X(s) := \int_{-\infty}^{\infty} e^{-st} X(t) dt
\eeq
where $\mathcal R\subset\C$ is the region of convergence of the
  integral.
For example, if $X(t)=e^{-3t}$, for $t\geq0$ and $X(t)=0$ otherwise, 
then $\hat X(s)=1/(s+3)$, and  
$\mathcal R = \{s=s_{re}+\jmath s_{im}:\ s_{re}>-3\}$
($\jmath$ is the imaginary number $\sqrt{-1}$). 

Some of its properties are:
\benu
\item For any constant matrix $A\in\R^{n\times n}$
\beq
   \widehat{AX}(s) = A\;\hat X(s);
\eeq
\item The Laplace transform of the derivative of $X$ is
\beq
    \widehat{\frac{dX}{dt}} (s) 
       = X(0)+s\;\int_{-\infty}^{\infty} e^{-st} X(t) dt 
       = X(0)+s\;\hat X (s)\ ;
\eeq
\item If $X(t+\delta)=:W(t)$ is a translation of $X$, then
\beq
    \hat W(s) = e^{-s\delta} \hat X(s)\ ;
\eeq
\item The inverse Laplace transform is 
\beq
    X(t) = \frac{1}{2\pi\jmath}\ 
          \int_{s_{re}-\jmath\infty}^{s_{re}+\jmath\infty} e^{st} \hat X(s) ds\, 
\eeq
with $s=s_{re}+\jmath s_{im}$, where $s_{re}$ is chosen so that $s_{re}+\jmath s_{im}$
is in the region of convergence $\mathcal R$.
\eenu

\subsubsection*{Transfer function}
Let $A\in\R^{n\times n}$, $B\in\R^{n\times m}$ and 
$C\in\R^{p\times n}$ be matrices, and let $X\in\R^n$, 
$Y\in\R^m$,  $R\in\R^p$, and consider the $n$-dimensional linear 
system with $m$ inputs and $p$ outputs:
\beqn{in-out-lin}
    \frac{dX}{dt} &=& AX +B R,\\
                   \nonumber \\
                Y &=& C\;X  \label{in-out-lin2}.
\eeqn
Applying the Laplace transform operator on both sides of 
the linear system\rf{in-out-lin}-(\ref{in-out-lin2}) yields an {\it
  algebraic equation} relating the new functions 
$\hat X(s)$, $\hat Y(s)$ and $\hat R(s)$:
\beq
    s \hat X(s) &=& A \hat X(s) + B\; \hat R(s)  \\
      \hat Y(s) &=& C\; \hat X(s).
\eeq
Moreover, for every $s$ for which the matrix $sI-A$ is invertible ($I$ is the
identity matrix),
\beq  
    (sI-A)\;\hat X(s) = B\; \hat R(s) \ \ \Rightarrow\ \ 
    \hat X(s) = (sI-A)^{-1}B\;\hat R(s)
\eeq
and thus, one can solve immediately for the output 
\beqn{out-fd}
    \hat Y(s) = C (sI-A)^{-1}B\;\hat R(s).
\eeqn
The transfer function of the system\rf{in-out-lin} is 
\beq
    \hat G(s) := C (sI-A)^{-1} B,
\eeq
and {\it depends only on the internal structure of the
system (i.e., $A$, $B$ and $C$)}.

\subsubsection*{Impulse response}
A useful case is that of the impulse response:
\beq
   R(t)=\delta(t),\ \ \ \ \ 
   \Rightarrow\ \ \hat R(s)\equiv1
\eeq
and therefore:
\beq
   \hat Y(s)\equiv \hat G(s)\ \ 
   \Leftrightarrow\ \ Y(t)\equiv  G(t),
\eeq
so that the transfer function of the system is the output corresponding to
a single pulse of input.

\subsubsection*{The gain $\dnorm{\hat G}_{\infty}$}
We have
\beq
    \dnorm{\hat Y}_2^2 = \frac{1}{2\pi}
        \int_{-\infty}^{\infty} 
                        \abs{\hat G(\jmath\omega) \hat R(\jmath\omega)}^2 d\omega
        \leq \frac{1}{2\pi}\sup_{\omega} \abs{\hat G(\jmath\omega)}^2 \ 
             \int_{-\infty}^{\infty}\abs{\hat R(\jmath\omega)}^2 d\omega 
\eeq
which is equivalent to
\beq
    \dnorm{\hat Y}_2 \leq \dnorm{\hat G}_{\infty} \dnorm{\hat R}_2
   \ \ \Leftrightarrow\ \ 
    \dnorm{Y}_2 \leq \dnorm{\hat G}_{\infty} \dnorm{R}_2.
\eeq
So, the infinity norm of the transfer function is an upper bound on
the strength of the output. 

To see that it is indeed the least upper bound, see for instance~\cite{control2}:
we can always choose a frequency $\omega_0$ so that 
\beq
  \dnorm{\hat G}_{\infty}=\abs{\hat G(\jmath\omega_0)}.
\eeq 
In our case, this is $\omega_0=0$. Then choose a control such that 
\beqn{eq-lt-sync}
  \abs{\hat R(\jmath\omega)}=\left\{\begin{array}{ll}
           r, & \mbox{ if } \abs{\omega}<\eps \\
           0, & \mbox{ otherwise, } 
         \end{array}\right.
\eeqn
where $\eps>0$ and $r$ should be such that $\hat R$ has unit 2-norm, 
for instance $r=\sqrt{\pi/\eps}$.
For very small $\eps>0$, $\abs{\hat R(\jmath\omega)}$ is zero, except on a very small 
neighborhood of $\omega_0=0$ and we may approximate:
\beq
   \frac{1}{2\pi} \int_{-\infty}^{\infty} \abs{\hat G(\jmath\omega)}^2
                            \abs{\hat R(\jmath\omega)}^2\ d\omega
&\approx& \frac{1}{2\pi}
          \int_{-\eps}^{\eps} r^2\ \abs{\hat G(\jmath\omega_0)}^2\ d\omega \\
& = & \frac{1}{2\pi}
      \abs{\hat G(\jmath\omega_0)}^2\ \int_{-\eps}^{\eps} r^2\ d\omega \\
& = & \dnorm{\hat G}_{\infty}^2 
\eeq
where the last equality follows from the definitions of $\omega_0$ and $r$.
Therefore
\beq
   \dnorm{Y}_2=\left[\frac{1}{2\pi}
                     \int_{-\infty}^{\infty} \abs{\hat G(\jmath\omega)}^2
                            \abs{\hat R(\jmath\omega)}^2\ d\omega 
               \right]^{\frac{1}{2}}
              \approx \dnorm{\hat G}_{\infty}.
\eeq
As an example of an input that (approximately)
satisfies\rf{eq-lt-sync}, consider $R(t)=2\frac{r}{\pi t}\sin{\eps t}$
(for $t\geq0$), the input plotted in Figure~\ref{fig-sync}.
Computation of the Laplace transform yields
$\hat R(s)=\frac{r}{\pi}\left[\pi-\mbox{Arctan}\frac{s}{\eps}
\right]$, where the function Arctan is the principal branch of the
complex inverse tangent function. It can be shown that, for
sufficiently small $\eps$, the function $\hat R$ approximately
satisfies condition\rf{eq-lt-sync}, except at the discontinuity points
$\omega=\pm\eps$.

\vspace{3mm}
More generally, in a system with $m$ inputs and $p$ outputs, one
defines the {\it internal gain of the system} $\dnorm{\hat G}_{\infty}$ 
as the induced $\mathcal L^2$ operator norm of the map from the inputs
to the outputs. It is possible to prove that
\beq
  \dnorm{\hat G}_{\infty} = \sup_{\omega\in\R} \bar\theta[\hat G(\jmath\omega)]
\eeq
where $\bar\theta$ denotes the largest singular value of the matrix 
$\hat G(\jmath\omega)$.

\subsubsection*{Stability of the transfer function}
As remarked above, expression\rf{out-fd} is valid if and only if  
the matrix $sI-A$ is invertible, or equivalently
\beq
     s \neq \lambda,\ \ \ \mbox{ for every eigenvalue, $\lambda$, of $A$}.
\eeq
If 
$\lambda_m=\max\{ \mbox{Re}(\lambda):\,
                 \lambda \mbox{ is an eigenvalue of $A$}\}$,
then the region of definition of the transfer function is included in the set  
$\mathcal R = \{s=s_{re}+\jmath s_{im}:\ s_{re}>\lambda_m \}$.

If all the eigenvalues of the matrix $A$ have negative real parts, then the transfer 
function is said to be stable. This is case for the matrix of the signaling 
cascade\rf{eq-in-out-sys}, whose eigenvalues are:
$
    -\beta_1,\ldots, -\beta_n,
$
so the transfer function $\hat G(s)$ is stable and well defined on 
$\mathcal R = \{s=s_{re}+\jmath s_{im}:\ s_{re}>-\min{\beta_i} \}$.

\end{document}